\documentclass[11pt,reqno]{amsart}
 
\usepackage[T1]{fontenc}
\usepackage[utf8]{inputenc}
\usepackage{lmodern}
\usepackage{amsmath,amssymb,amsthm,mathtools}
\usepackage{enumitem}
\setlist[enumerate]{font=\upshape}
\setlist[enumerate,2]{leftmargin=2.5em,itemsep=0.3em,topsep=0.3em}
\usepackage{tikz}
\usepackage{caption}
\usepackage[a4paper,left=3cm,right=3cm,top=3cm,bottom=3cm]{geometry}
\usepackage{hyperref}

\setlist[enumerate,1]{label=\textup{(\roman*)}}
\newtheorem{theorem}{Theorem}[section]
\newtheorem{proposition}[theorem]{Proposition}
\newtheorem{lemma}[theorem]{Lemma}
\newtheorem{corollary}[theorem]{Corollary}
\theoremstyle{definition}
\newtheorem{definition}[theorem]{Definition}
\theoremstyle{remark}
\newtheorem{remark}[theorem]{Remark}
\numberwithin{equation}{section}

\begin{document}

\title[Spectral configurations under rank-one perturbations]%
{Prescribed spectral configurations under rank-one perturbations}

\author[M. Hansmann]{Marcel Hansmann}
\address{Faculty of Mathematics, Chemnitz University of Technology,
09107 Chemnitz, Germany}
\email{marcel.hansmann@mathematik.tu-chemnitz.de}

\begin{abstract}
  A spectral configuration specifies the Fredholm set, the Fredholm index on each of its components, and, on zero-index components, either the whole component or a locally finite set of discrete eigenvalues with prescribed algebraic multiplicities. We show that any two spectral configurations with the same Fredholm set and index can be realized by a suitable closed operator \(A\) and a rank-one perturbation \(A+K\). This remains true even when the rank-one operator \(K\), acting on a given separable infinite-dimensional Hilbert space, is prescribed in advance. In particular, for any such prescribed \(K\), one can choose \(A\) such that the discrete eigenvalues of \(A+K\) accumulate at the common essential spectrum at any prescribed rate.

The proof is inspired by classical linearization theory. The main tool is the fact that arbitrary holomorphic functions can be represented as rank-one perturbation determinants.
\end{abstract}

 \maketitle

\section{Introduction and main results}\label{sec:intro}

Given two arbitrary families \(\lambda_1,\ldots,\lambda_n\) and \(\mu_1,\ldots,\mu_n\) of complex numbers, it is easy to construct a matrix \(A_n\in\mathbb C^{n\times n}\) and a matrix \(K_n\) of rank at most one such that the spectra of \(A_n\) and \(A_n+K_n\) consist precisely of these numbers, counting algebraic multiplicities. Indeed, one may simply take the companion matrices associated with the polynomials \(p(z)=\prod_{j=1}^{n}(z-\lambda_j)\) and \(q(z)=\prod_{j=1}^{n}(z-\mu_j)\), respectively; see Krupnik~\cite{Krupnik1992}. These matrices can differ only in one row and hence by a matrix of rank at most one.

Thus, in finite dimensions, by choosing \(A_n\) suitably, a perturbation of rank at most one suffices to transform an arbitrary spectral configuration into any other one. In this article we establish an infinite-dimensional counterpart of this result.
\medskip

\noindent\textbf{Spectral configurations and the main result.}
To explain what we are aiming for, we first need to discuss what we mean by a spectral configuration in infinite dimensions.

Let us start with a densely defined closed operator \(A\) on a separable Hilbert space \(\mathcal H\). We denote its spectrum and resolvent set by \(\sigma(A)\) and \(\rho(A)\), respectively.
Moreover,
\[
  \Phi(A):=\{\lambda\in\mathbb C:\lambda-A\text{ is a Fredholm operator}\}
\]
denotes the Fredholm set of \(A\), and
\(\sigma_{\mathrm{ess}}(A):=\mathbb C\setminus\Phi(A)\) denotes its
Fredholm essential spectrum. Recall that \(\lambda-A\) is called Fredholm
if its kernel is finite dimensional and its range is closed and of
finite codimension. The Fredholm set is open and the index function
\[
  \lambda\mapsto
  \operatorname{ind}(\lambda-A)
  :=
  \dim\ker(\lambda-A)
  -
  \operatorname{codim}\operatorname{Ran}(\lambda-A)
\]
is locally constant on \(\Phi(A)\). In particular, on each connected
component \(U\) of \(\Phi(A)\), it has a constant value \(\kappa_U\),
which constrains the possible spectral behavior of \(A\) in \(U\).

\begin{enumerate}
\item
If \(\kappa_U\neq0\), then \(\lambda-A\) cannot be invertible for any \(\lambda\in U\), and hence \(U\subseteq\sigma(A)\), i.e. the entire component is contained in the spectrum.

\item If \(\kappa_U=0\), there are two possibilities: either
\(U\subseteq\sigma(A)\), or \(\sigma(A)\cap U\) is locally finite
in \(U\), that is, it intersects every compact subset of \(U\)
in at most finitely many points. In the latter case every
\(\lambda\in\sigma(A)\cap U\) is a \emph{discrete eigenvalue}
of \(A\), i.e., an isolated eigenvalue of finite algebraic
multiplicity \(m_a(\lambda;A)\), see Definition~\ref{def:riesz-projection}. If there are infinitely
many such eigenvalues, they can have finite accumulation points
only on \(\partial U\).
\end{enumerate}

Most important for us, the Fredholm properties of a closed operator are not affected by compact perturbations. More precisely, for a compact operator \(K\) on \(\mathcal{H}\) and \(\lambda\in\mathbb C\), \(\lambda-A\) is Fredholm if and only if \(\lambda-(A+K)\) is Fredholm, and in this case their Fredholm indices coincide. Consequently, \(\sigma_{\mathrm{ess}}(A)=\sigma_{\mathrm{ess}}(A+K)\). For all these standard facts from Fredholm theory, we refer to Kato~\cite[Chapter~IV, Section~5]{Kato} or Gohberg et al.~\cite[Chapter~XVII]{GohbergGoldbergKaashoek1990}.

\begin{definition}\label{def:spectral-configuration}
Let \(\Omega\subseteq\mathbb C\) be open.
\begin{enumerate}
\item
A \emph{spectral configuration} on \(\Omega\) is a tuple \((\Omega_{\mathrm{full}},\kappa,\Lambda,m)\) such that
\begin{enumerate}
\item \(\Omega_{\mathrm{full}}\) is a possibly empty union of connected components of \(\Omega\).

\item \(\kappa:\Omega\to\mathbb Z\) is locally constant and \(\{\lambda\in\Omega:\kappa(\lambda)\neq0\} \subseteq\Omega_{\mathrm{full}}\).

\item \(\Lambda\subseteq\Omega\setminus\Omega_{\mathrm{full}}\) is locally finite in \(\Omega\).

\item \(m:\Lambda\longrightarrow\mathbb N\).
\end{enumerate}

\item A closed operator \(A\) \emph{realizes} the spectral configuration \((\Omega_{\mathrm{full}},\kappa,\Lambda,m)\) \emph{on \(\Omega\)} if
\begin{enumerate}
\item \(\Phi(A)=\Omega\) or, equivalently, \(\sigma_{\mathrm{ess}}(A)=\mathbb C\setminus\Omega\).

\item \(\sigma(A)\cap\Omega=\Omega_{\mathrm{full}}\cup\Lambda\).

\item \(\operatorname{ind}(\lambda-A)=\kappa(\lambda)\) for all \(\lambda\in\Omega\).

\item Every \(\lambda\in\Lambda\) is a discrete eigenvalue of \(A\) with \(m_a(\lambda;A)=m(\lambda)\).
\end{enumerate}
\end{enumerate}
\end{definition}

Let us briefly explain the meaning of the four entries of a spectral configuration. The set \(\Omega_{\mathrm{full}}\) consists of those components of \(\Omega=\Phi(A)=\mathbb C\setminus\sigma_{\mathrm{ess}}(A)\) which are prescribed to lie entirely in the spectrum, while \(\kappa\) prescribes the Fredholm index. Thus every component on which \(\kappa\) is nonzero necessarily belongs to \(\Omega_{\mathrm{full}}\). On the components outside \(\Omega_{\mathrm{full}}\), the prescribed spectrum is the set \(\Lambda\), which is locally finite in \(\Omega\), and \(m(\lambda)\) prescribes the algebraic multiplicity of the corresponding discrete eigenvalue. Notice that \(\Omega_{\mathrm{full}}\) may also contain components on which \(\kappa\) vanishes. Thus, on a zero-index component, the prescribed spectrum may be the whole component or an arbitrary locally finite set of discrete eigenvalues, including the empty set.

\begin{remark}\label{rem:spectral-picture}
For bounded operators, Definition~\ref{def:spectral-configuration} is closely related to the classical notion of a \emph{spectral picture}, terminology introduced by Pearcy~\cite{Pearcy1978}. See also Chevreau~\cite{Chevreau1980}. In Chevreau's formulation, the spectral picture consists of the set on which \(\lambda-A\) fails to be semi-Fredholm together with the extended semi-Fredholm index on its complement, taking values in \(\mathbb Z\cup\{-\infty,+\infty\}\). Thus it distinguishes the semi-Fredholm regions, including those of infinite index, but does not record the finer spectral disposition on Fredholm components of index zero, such as whether an entire component is contained in the spectrum or which discrete eigenvalues, with which algebraic multiplicities, occur in that component. Conway~\cite{Conway1985} subsequently studied precisely this finer problem, again for a single bounded operator.

Our spectral configurations are closer in spirit to Conway's viewpoint. As far as the semi-Fredholm structure is concerned, however, we retain only the Fredholm/non-Fredholm decomposition and the ordinary Fredholm index on the Fredholm region. Since our main problem is the simultaneous realization of \emph{two} spectral configurations rather than one, we avoid the additional layer of complexity that would arise from prescribing semi-Fredholm data and infinite indices as well.
\end{remark}

\begin{figure}[ht]
  \centering
  \begin{tikzpicture}[scale=0.6]
    \begin{scope}[xshift=-5.2cm]
      \fill[gray!20] (-2,-1) circle (1.2);
      \fill[gray!20] (2,-1) circle (1);

      \draw[thick] (-2,-1) circle (1.2);
      \draw[thick] (2,-1) circle (1);
      \draw[thick] (0,2) circle (1.5);

      \draw[red,line width=0.8pt] (0.28,1.88) -- (0.42,2.02);
      \draw[red,line width=0.8pt] (0.28,2.02) -- (0.42,1.88);

      \draw[red,line width=0.8pt] (0.52,2.18) -- (0.66,2.32);
      \draw[red,line width=0.8pt] (0.52,2.32) -- (0.66,2.18);

      \draw[red,line width=0.8pt] (-0.18,2.38) -- (-0.04,2.52);
      \draw[red,line width=0.8pt] (-0.18,2.52) -- (-0.04,2.38);

      \draw[red,line width=0.8pt] (-0.52,1.62) -- (-0.38,1.76);
      \draw[red,line width=0.8pt] (-0.52,1.76) -- (-0.38,1.62);

      \draw[red,line width=0.8pt] (0.02,1.42) -- (0.16,1.56);
      \draw[red,line width=0.8pt] (0.02,1.56) -- (0.16,1.42);

      \draw[red,line width=0.8pt] (2.88,2.28) -- (3.02,2.42);
      \draw[red,line width=0.8pt] (2.88,2.42) -- (3.02,2.28);

      \draw[red,line width=0.8pt] (3.08,1.48) -- (3.22,1.62);
      \draw[red,line width=0.8pt] (3.08,1.62) -- (3.22,1.48);

      \draw[red,line width=0.8pt] (-3.02,1.18) -- (-2.88,1.32);
      \draw[red,line width=0.8pt] (-3.02,1.32) -- (-2.88,1.18);

      \draw[red,line width=0.8pt] (2.68,-2.02) -- (2.82,-1.88);
      \draw[red,line width=0.8pt] (2.68,-1.88) -- (2.82,-2.02);

      \node at (-2,-1.25) {\(\kappa=1\)};
      \node at (2,-1.25) {\(\kappa=0\)};
      \node at (0,1.25) {\(\kappa=0\)};
      \node at (-2.85,2.65) {\(\kappa=0\)};

      \node (Afull) at (0,-2.6) {\(\Omega_{A,\mathrm{full}}\)};
      \draw[->]
        ([xshift=-4mm]Afull.north)
        --
        (-1.35,-1.55);
      \draw[->]
        ([xshift=4mm]Afull.north)
        --
        (1.4,-1.55);

      \node at (0,4.0) {\(\sigma(A)\)};
    \end{scope}

    \draw[->,very thick] (-0.9,0.7) -- (0.9,0.7);
    \node at (0,1.05) {\(+\,K\)?};

    \begin{scope}[xshift=5.2cm]
      \fill[gray!20] (-2,-1) circle (1.2);

      \fill[gray!20] (0,2) circle (1.5);

      \draw[thick] (-2,-1) circle (1.2);
      \draw[thick] (2,-1) circle (1);
      \draw[thick] (0,2) circle (1.5);

      \draw[blue,line width=0.8pt] (1.55,-0.85) -- (1.69,-0.71);
      \draw[blue,line width=0.8pt] (1.55,-0.71) -- (1.69,-0.85);


      \draw[blue,line width=0.8pt] (2.38,-0.72) -- (2.52,-0.58);
      \draw[blue,line width=0.8pt] (2.38,-0.58) -- (2.52,-0.72);

      \draw[blue,line width=0.8pt] (1.82,-1.05) -- (1.96,-0.91);
      \draw[blue,line width=0.8pt] (1.82,-0.91) -- (1.96,-1.05);

      \draw[blue,line width=0.8pt] (3.02,2.58) -- (3.16,2.72);
      \draw[blue,line width=0.8pt] (3.02,2.72) -- (3.16,2.58);

      \draw[blue,line width=0.8pt] (3.20,0.55) -- (3.34,0.69);
      \draw[blue,line width=0.8pt] (3.20,0.69) -- (3.34,0.55);

      \draw[blue,line width=0.8pt] (-3.15,1.52) -- (-3.01,1.66);
      \draw[blue,line width=0.8pt] (-3.15,1.66) -- (-3.01,1.52);

      \draw[blue,line width=0.8pt] (0.85,-1.95) -- (0.99,-1.81);
      \draw[blue,line width=0.8pt] (0.85,-1.81) -- (0.99,-1.95);

      \node at (-2,-1.25) {\(\kappa=1\)};
      \node at (2,-1.25) {\(\kappa=0\)};
      \node at (0,1.25) {\(\kappa=0\)};
      \node at (-2.85,2.65) {\(\kappa=0\)};

      \node (Bfull) at (0,-2.8) {\(\Omega_{A+K,\mathrm{full}}\)};
      \draw[->]
        ([xshift=-5mm]Bfull.north)
        --
        (-1.45,-1.55);
      \draw[->]
        ([xshift=5mm]Bfull.north)
        --
        (0.55,0.85);

      \node at (0,4.0) {\(\sigma(A+K)\)};
    \end{scope}
  \end{tikzpicture}

  \caption{Two spectral configurations. The circles form the essential
  spectrum, the shaded discs are components contained entirely in the
  spectrum, and the crosses mark discrete eigenvalues.}
  \label{fig:spectral-configurations}
\end{figure}
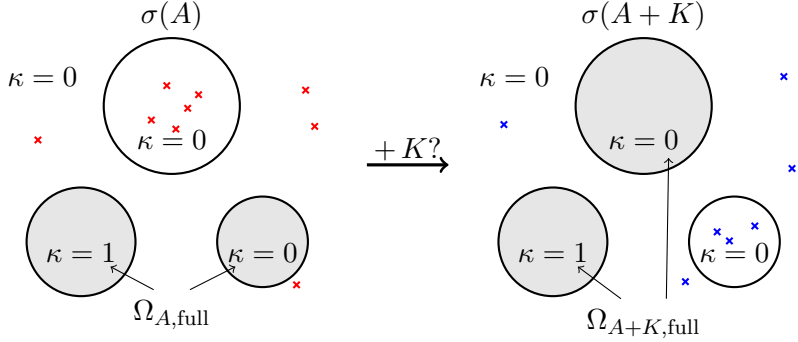

With this terminology in place, we now ask: how can the spectral configurations of \(A\) and \(A+K\) differ if \(K\) has rank one? For instance, can the configurations shown in Figure~\ref{fig:spectral-configurations} be realized by an actual pair of operators?

\medskip

Note that in Figure~\ref{fig:spectral-configurations} the three circles, representing the common essential spectrum, are the same for \(A\) and \(A+K\), and so are the values of the Fredholm index on the components of the Fredholm set. These are exactly the obvious necessary restrictions imposed by Fredholm stability. Are these restrictions also sufficient? More precisely, could the spectra of \(A\) and \(A+K\) on zero-index components be prescribed independently, with the spectrum of each operator on such a component being either the entire component or an arbitrary locally finite set of discrete eigenvalues, with arbitrary algebraic multiplicities?

Our main result shows that this is indeed the case, and it remains true even when the Hilbert space and the rank-one perturbation are prescribed in advance.

\begin{theorem}\label{thm:main}
Let \(\Omega\subseteq\mathbb C\) be open and let
\[
  \mathcal C_A=(\Omega_{A,\mathrm{full}},\kappa,\Lambda_A,m_A)
  \quad\text{and}\quad
  \mathcal C_B=(\Omega_{B,\mathrm{full}},\kappa,\Lambda_B,m_B)
\]
be two spectral configurations on \(\Omega\) with the same index function \(\kappa:\Omega\to\mathbb Z\). Moreover, let \(\mathcal H\) be a separable infinite-dimensional Hilbert space and let \(K\in\mathcal B(\mathcal H)\) have rank one. Then there exists a densely defined closed operator \(A\) on \(\mathcal H\) such that \(A\) realizes \(\mathcal C_A\) on \(\Omega\) and \(B=A+K\) realizes \(\mathcal C_B\) on \(\Omega\).
\end{theorem}

\begin{remark}
  \label{rem:2}
Let us note that prescribing the Hilbert space in Theorem~\ref{thm:main} would by itself be immaterial, since all separable infinite-dimensional complex Hilbert spaces are unitarily isomorphic. The additional point is that the concrete rank-one perturbation \(K\) may be prescribed. This will be achieved by first constructing a rank-one operator \(K'\) with the same trace as \(K\) and then passing from \(K'\) to \(K\) by a bounded similarity, see Lemma~\ref{lem:similarity}.
\end{remark}

Two spectral configurations can be realized by bounded operators precisely when the spectra are bounded and the common essential spectrum is nonempty. Moreover, one can also prescribe perturbations of general finite rank.

\begin{corollary}\label{cor:bounded-main}
In the situation of Theorem~\ref{thm:main}, the operator \(A\) can be
chosen in \(\mathcal B(\mathcal H)\) (and then
\(B=A+K\in\mathcal B(\mathcal H)\) as well) if and only if
\(F:=\mathbb C\setminus\Omega\) is nonempty and the set
\[
  F\cup\Omega_{A,\mathrm{full}}\cup\Lambda_A
  \cup\Omega_{B,\mathrm{full}}\cup\Lambda_B
\]
is bounded.
\end{corollary}

\begin{corollary}\label{cor:finite-rank}
The assertions of Theorem~\ref{thm:main} and Corollary~\ref{cor:bounded-main} remain valid for every prescribed \(K\in\mathcal B(\mathcal H)\) with \(0<\operatorname{rank}K<\infty\).
\end{corollary}

\medskip

\noindent\textbf{Scope and limitations.}
Let us illustrate the scope of these results by a few examples. On \(\ell^2(\mathbb N)\), fix \(\alpha\in\mathbb C\setminus\{0\}\) and let
\[
  K=
  \begin{pmatrix}
    \alpha & 0 &  \cdots\\
    0      & 0 &  \cdots\\
    \vdots & \vdots &  \ddots
  \end{pmatrix}.
\]
That is, \(Kf=\alpha\langle f,e_1\rangle e_1\), where \(e_1=(1,0,\ldots)\). Then the following holds:
\begin{enumerate}
\item There exists \(A_1\) such that \(A_1\) and \(A_1+K\) realize the two spectral configurations indicated in Figure~\ref{fig:spectral-configurations}. \item As a more complicated example, there exists \(A_2\) such that
\[
  \sigma_{\mathrm{ess}}(A_2)
  =
  \sigma_{\mathrm{ess}}(A_2+K)
  =
  \bigcup_{n\in\mathbb Z}(n+i\mathbb R),
\]
with Fredholm index zero on every complementary strip. On each such strip, independently for \(A_2\) and \(A_2+K\), the spectrum may be prescribed to be either the entire component or an arbitrary locally finite set of discrete eigenvalues with arbitrary algebraic multiplicities. \item At the opposite extreme, there is also an operator \(A_3\) such that
\[
  \sigma(A_3)
  =
  \sigma_{\mathrm{ess}}(A_3)
  =
  \sigma_{\mathrm{ess}}(A_3+K)
  =
  \emptyset,
  \qquad
  \sigma(A_3+K)=\mathbb C.
\]
\end{enumerate}

On the other hand, Theorem~\ref{thm:main} also has natural limitations. For instance, it does not prescribe arbitrary Jordan structures. In fact, all prescribed discrete eigenvalues of the realizations constructed below have geometric multiplicity one. For \(\lambda\in\Lambda_A\cap\rho(B)\) (resp. \(\lambda\in\Lambda_B\cap\rho(A)\)), this is in fact necessary. To see this, fix \(\lambda\in\Phi(A)\) and let \(E\) be the intersection of \(\ker(\lambda-A)\) and \(\ker(\lambda-A-K)\). Then \(E\) is the kernel of each of the restrictions \(K|_{\ker(\lambda-A)}\) and \(K|_{\ker(\lambda-A-K)}\). Using the rank-nullity theorem and the fact that both of these restrictions have rank at most one, we obtain
\[
\bigl|\dim\ker(\lambda-A)-\dim\ker(\lambda-A-K)\bigr|\le1.
\]
This elementary obstruction is the rank-one instance of more general bounds on root subspaces and Jordan chains under finite-rank perturbations; see Behrndt et al.~\cite{BehrndtLebenMartinezPeriaTrunk2015}. Beyond such necessary restrictions, we do not know to what extent Jordan structures can be prescribed simultaneously in addition to the spectral configurations. For comparison, in finite dimensions Thompson~\cite{Thompson1980} completely describes, for a fixed matrix \(A\), the invariant factors, and hence the Jordan structures, attainable by a perturbation of \(A\) of rank at most one.

Another limitation concerns the admissible perturbations. Corollary~\ref{cor:finite-rank} covers every prescribed nonzero finite-rank perturbation, but we do not know whether this extends to arbitrary prescribed compact operators of infinite rank as well.

\medskip
\noindent \textbf{Consequences for eigenvalue estimates.} In Corollary~\ref{cor:rate} we show that, for any prescribed
 rank-one perturbation \(K\) and any nonempty proper closed set \(F\subseteq\mathbb C\), one may choose a closed operator \(A\) with \(\sigma(A)=\sigma_{\mathrm{ess}}(A)=F\) such that, for every \(p>0\),
\begin{equation}\label{eq:1}
\sum_{\lambda\in\sigma_d(A+K)} \operatorname{dist}(\lambda,\sigma(A))^p=\infty.
\end{equation}
If \(F\) is compact, \(A\) may be chosen bounded. Consequently, a finite-rank assumption alone cannot yield general Lieb--Thirring-type summability estimates, or comparable Blaschke-type restrictions, relative to the (essential) spectrum. This complements the many positive results of this kind and explains why these require additional information on the unperturbed operator, such as selfadjointness, normality or resolvent estimates. For abstract eigenvalue variation estimates under selfadjointness or normality assumptions, see, e.g., \cite{Kato1987,BhatiaElsner1994,BhatiaDavis1999,HansmannVariation2013}; for related Lieb--Thirring- and Blaschke-type bounds in non-selfadjoint settings, see, for instance, \cite{BruneauOuhabaz2008,BorichevGolinskiiKupin2009, FavorovGolinskii2009,DemuthHansmannKatriel2009,FavorovGolinskii2015, GolinskiiKupin2017,Frank2018, BoegliPetpradittha2026}. Let us note, however, that in the bounded case the sum in \eqref{eq:1} is finite for every \(p\geq1\) if \(\sigma(A)\) is replaced by the numerical range of \(A\), see \cite{HansmannNumericalRange2011}.

\medskip

\noindent\textbf{A scalar resolvent realization and the idea of proof.}
One of the main ingredients in the proof of Theorem~\ref{thm:main} is the following realization result for holomorphic functions.

\begin{proposition}
\label{prop:scalar-realization}
Let \(\Omega\subseteq\mathbb C\) be open and nonempty and let \(g\in H(\Omega)\). Then there exist a separable infinite-dimensional Hilbert space \(\mathcal G\), a densely defined closed operator \(T\) on \(\mathcal G\), and vectors \(x,y\in\mathcal G\), with \(x\ne0\) and \(\|y\|=1\), such that
\begin{equation}\label{eq:3}
  \rho(T)=\Omega,
  \qquad
  \sigma(T)=\sigma_{\mathrm{ess}}(T)=\mathbb C\setminus\Omega,
\end{equation}
and
\[
  g(\lambda)
  =
  \langle(\lambda-T)^{-1}x,y\rangle,
  \qquad \lambda\in\Omega.
\]
\end{proposition}

In fact, Theorem~\ref{thm:simultaneous} below gives a stronger version of Proposition~\ref{prop:scalar-realization}, in which finitely many holomorphic functions \(g_1,\ldots,g_n\) can be realized simultaneously with a common operator \(T\), a common vector \(y\), and linearly independent representing vectors \(x_1,\ldots,x_n\). Moreover, the pairings \(\langle x_i,y\rangle\) can be prescribed in advance.

The relevance of this proposition to our problem is immediate from the rank-one determinant identity. Indeed, given \(h\in H(\Omega)\), apply the proposition to \(g=1-h\) and set \(Kf=\langle f,y\rangle x\), \(f\in\mathcal G\). Then
\[
h(\lambda)=\det\bigl(I-K(\lambda-T)^{-1}\bigr), \qquad \lambda\in\Omega=\rho(T),
\]
i.e. \(h\) is the perturbation determinant of \(T\) by \(K\). In particular, \(h(\lambda)=0\) if and only if \(\lambda\in\sigma(T+K)\). Moreover, every isolated zero of \(h\) corresponds to a discrete eigenvalue of \(T+K\), with its order equal to the algebraic multiplicity. Thus, to realize the spectral data on the zero-index components of a given spectral configuration, we can take \(h\equiv0\) on full components and use Weierstrass' theorem to prescribe its zeros on the remaining ones.

For the simultaneous realization of two spectral configurations, we then choose two such holomorphic functions and use Theorem~\ref{thm:simultaneous} to realize them with a common operator \(T\) and, crucially, a common vector \(y\). This produces two rank-one perturbations
\[
K_A=\langle\,\cdot\,,y\rangle x_A,
\qquad
K_B=\langle\,\cdot\,,y\rangle x_B,
\]
so that
\[
(T+K_B)-(T+K_A)
=\langle\,\cdot\,,y\rangle(x_B-x_A)
\]
has rank one, since \(x_A\) and \(x_B\) are linearly independent. We
prescribe the pairings \(\langle x_A,y\rangle\) and
\(\langle x_B,y\rangle\) so that this difference has the same trace as
the prescribed perturbation \(K\).

A separate operator then carries the prescribed Fredholm index on the nonzero-index components, see Proposition~\ref{prop:index}. Taking a common direct sum combines the zero-index realization with this index carrier. Finally, Lemma~\ref{lem:similarity} transfers the resulting rank-one difference to the prescribed operator \(K\).

\begin{remark}
\label{rem:1}
To the best of our knowledge, Proposition~\ref{prop:scalar-realization} in this exact form does not appear in the literature. Nevertheless, we do not regard the underlying realization mechanism as new. For bounded open sets it is contained, at least implicitly, in classical linearization theory. Mitiagin's construction has explicitly computable scalar resolvent
coefficients \cite[Lemma~3.2]{Mityagin1978II} and, with a common
choice of weight, also allows finitely many holomorphic functions
to be realized simultaneously; see also
\cite[Proposition~4.8]{Mityagin1978II} for the Hilbert-space setting.
Den Boer proves a closely related evaluation identity
\cite[Lemma~3(iii)]{DenBoer1978}. See also Gohberg et al.~\cite{GohbergKaashoekLay1978}. For arbitrary open sets, den Boer works with linear pencils \(S-\lambda V\) rather than with the resolvent of a single closed operator.

In Proposition~\ref{prop:scalar-realization} we use a formulation adapted to rank-one perturbation determinants together with the exact spectral normalization \eqref{eq:3}. We do not use the classical results in the proof. Instead, we give a direct, self-contained construction along similar lines, which also provides the Fredholm-index carriers needed later.
\end{remark}

Let us mention that scalar resolvent coefficients as in Proposition~\ref{prop:scalar-realization} are also known as \emph{weak resolvents}, see \cite{NordgrenRadjaviRosenthal1987,FongNordgrenRadjaviRosenthal1990}. Fixed resolvent coefficients together with exact spectral normalizations also occur in the theory of spectral nodes, see \cite{KaashoekVanDerMeeRodman1981}.

In earlier joint work with Katriel~\cite{HansmannKatriel2013}, a related resolvent-coefficient representation was used for Cauchy transforms on the unit disk in order to relate their zeros to eigenvalues of rank-one perturbations.

\medskip

\noindent\textbf{Related spectral problems.}
For a single bounded operator, the realization problem behind
Definition~\ref{def:spectral-configuration} is classical, as discussed
in Remark~\ref{rem:spectral-picture}. What is not contained in these
results is a simultaneous realization of two spectral configurations
\(\mathcal{C}_A\) and \(\mathcal{C}_B\) in which the two operators
are coupled by the requirement that \(B-A\) be a prescribed rank-one
operator.

Finite-dimensional rank-one spectral perturbations, including
structured classes of matrices, are studied by
Krupnik~\cite{Krupnik1992}. In this case, if the perturbation \(K_n\) itself is
prescribed in addition to the two spectra, one has the necessary
trace identity
\[
  \operatorname{tr}(A_n+K_n)-\operatorname{tr}A_n
  =\operatorname{tr}K_n.
\]
By contrast, Theorem~\ref{thm:main} and
Corollary~\ref{cor:bounded-main} show that no corresponding global
trace compatibility condition arises in infinite dimensions.

Theorem~\ref{thm:main} also differs from classical spectral assignment
problems in what is kept fixed: there, the background operator is fixed
and one studies which spectra can be obtained by varying the
perturbation. In addition to the papers concerned with eigenvalue estimates mentioned above, see, e.g.,
Gurvits et al.~\cite{GurvitsRodmanSpitkovsky1991} and, for
rank-one perturbations of a fixed selfadjoint operator with discrete
spectrum,
Dobosevych--Hryniv~\cite{DobosevychHrynivLAA,DobosevychHrynivIEOT}.
In the selfadjoint setting, restrictions persist even when both
operators may vary: for cyclic selfadjoint operators with pure point
spectrum, Poltoratski~\cite{Poltoratski1998} obtains
nontrivial compatibility conditions on the two eigenvalue sets.
Once selfadjointness is dropped and only the spectral configurations
are prescribed, Theorem~\ref{thm:main} shows that Fredholm stability
is the only restriction.

\medskip

\noindent\textbf{Contents of the paper.}
In Section~\ref{sec:ingredients} we introduce a multiplication-operator model and prove Proposition~\ref{prop:scalar-realization} and its simultaneous generalization. In Section~\ref{sec:pair}, we construct an index carrier for the prescribed Fredholm indices and combine it with rank-one perturbation determinants to prove our main results. The direct-sum facts used in the proofs are collected in Appendix~\ref{app:direct-sums}. A \hyperlink{ai-declaration}{declaration} at the end of the paper describes the role of generative AI in this work.

\medskip

\noindent\textbf{A word on notation.}
Throughout this text, all Hilbert spaces are complex, and their inner products \(\langle\cdot,\cdot\rangle\) are taken to be linear in the first argument. We use the same symbol for the inner products on different Hilbert spaces, their meaning always being clear from the context. The same applies to the symbol \(\|\cdot\|\), which is used both for norms on various spaces and for operator norms. For Hilbert spaces \(\mathcal H_1,\mathcal H_2\), \(\mathcal B(\mathcal H_1,\mathcal H_2)\) denotes the space of bounded linear operators from \(\mathcal H_1\) to \(\mathcal H_2\), and we write \(\mathcal B(\mathcal H):=\mathcal B(\mathcal H,\mathcal H)\). For a linear operator \(A\), its domain is denoted by \(\mathcal D(A)\). Its Fredholm set \(\Phi(A)\) and essential spectrum \(\sigma_{\mathrm{ess}}(A)\) have already been defined above. We write \(\sigma_d(A)\) for its discrete spectrum. For a densely defined operator \(T:\mathcal D(T)\subseteq\mathcal H_1\to\mathcal H_2\), its Hilbert space adjoint is denoted by \(T^*:\mathcal D(T^*)\subseteq\mathcal H_2\to\mathcal H_1\).

We denote planar Lebesgue measure by \(dA\). For an open set \(G\subseteq\mathbb C\), \(H(G)\) denotes the space of holomorphic functions on \(G\). Moreover, \(G^\#:=\{\overline z:z\in G\}\) denotes the reflection of \(G\) across the real axis. When no confusion arises, functions are identified with their defining expressions. For instance, \(z,z^2,\ldots\) denote the functions \(w\mapsto w,w^2,\ldots\) on a given domain. By \(\mathbf 1\) we denote the constant function with value \(1\).

\section{The resolvent realization}\label{sec:ingredients}

In this section we prove Proposition~\ref{prop:scalar-realization} and its simultaneous generalization. As indicated in Remark~\ref{rem:1} and the references cited therein, our construction draws on and adapts classical ideas from linearization theory. Before turning to this construction, we indicate the main idea by a few formal computations.

We will use the following difference quotient operator. Note that the apparent singularity at \(z=\lambda\) is removable.

\begin{definition}
  \label{def:1}
Let \(G\subseteq\mathbb C\) be open and nonempty. For \(\lambda \in G\), define \(D_{\lambda}: H(G) \to H(G)\) by
\[
  (D_\lambda f)(z)
  :=
  \begin{cases}
    \dfrac{f(z)-f(\lambda)}{z-\lambda},
    & z\neq\lambda,\\[2mm]
    f'(\lambda),
    & z=\lambda.
  \end{cases}
\]
\end{definition}

Now suppose that \(M\) is the coordinate multiplication operator on a suitable Hilbert space \(\mathcal K\subseteq H(G)\), that is,
\[
    (Mf)(z)=zf(z), \qquad z\in G,\quad f\in\mathcal D(M).
  \]
Then, at least formally, for all \(\lambda \in G\) we have \((-D_{\lambda})(\lambda-M)=I\) and
\[
  (\lambda-M)(-D_\lambda f)=f-f(\lambda)\mathbf1.
\]
Thus \(\lambda-M\) is injective and has a one-dimensional cokernel, since its range consists precisely of the functions vanishing at \(\lambda\). If a second operator \(N\) on a Hilbert space \(\mathcal K_0\) has the
complementary defect for every \(\lambda\in G\), namely if
\(\lambda-N\) is surjective with a one-dimensional kernel, and
\(v_\lambda\in\ker(\lambda-N)\) is normalized by
\(\langle v_\lambda,e\rangle=1\) for a fixed unit vector
\(e\in\mathcal K_0\), then the rank-one upper triangular coupling
\[
  T=
  \begin{pmatrix}
    M&\langle\,\cdot\,,e\rangle\mathbf1\\
    0&N
  \end{pmatrix}
\]
cancels the two defects, hence \(\lambda-T\) will be invertible for
\(\lambda\in G\) and its inverse contains the value \(f(\lambda)\) as a
fixed scalar coefficient.  Indeed, at least formally, we have
\[
    (\lambda-T)
    (-D_\lambda f, -f(\lambda)v_\lambda)
    = (f,0)
\quad \text{ and so } \quad
    (\lambda-T)^{-1}(f,0)
    = (-D_\lambda f, -f(\lambda)v_\lambda).
\] 
Taking the inner product with the fixed vector \((0,-e)\) therefore gives
\[
    f(\lambda)
    =
    \left\langle
        (\lambda-T)^{-1}(f,0),(0,-e)\right\rangle.
  \]
  
In the remainder of this section we make these formal computations precise. We begin by constructing the Hilbert space \(\mathcal K\subseteq H(G)\) and the multiplication operator \(M\). The latter will also be used in the construction of the Fredholm-index carriers in Proposition~\ref{prop:index}. Later, \(N\) will be chosen as the adjoint of the corresponding multiplication operator on the reflected domain \(G^\#\).

\begin{lemma}\label{lem:multiplication}
Let \(G\subseteq\mathbb C\) be open and nonempty and let \(a_1,\ldots,a_s\in H(G)\), where \(s\in\mathbb N_0\). Then there is a separable infinite-dimensional Hilbert space
\[
  \mathcal K=\mathcal K(a_1,\ldots,a_s,G)\subseteq H(G)
\]
(simply denoted by \(\mathcal{K}(G)\) if \(s=0\)) with the following properties.
\begin{enumerate}
\item
Every polynomial and each \(a_j\) belong to \(\mathcal K\), and for every \(\lambda\in G\) both the point evaluation
\[
  \delta_{\lambda} : \mathcal{K} \to \mathbb{C}, \quad \delta_\lambda f=f(\lambda),
\]
and the difference quotient \(D_\lambda: \mathcal{K} \to \mathcal{K}\) are bounded on \(\mathcal{K}\).

\item
The maximal coordinate multiplication operator
\[
  (Mf)(z)=zf(z),
  \qquad
  \mathcal D(M)=\{f\in\mathcal K:zf\in\mathcal K\},
\]
is densely defined and closed on \(\mathcal{K}\).

\item
Let \(\lambda\in G\). Then
\[
  (-D_\lambda)(\lambda-M)=I
  \quad\text{on }\mathcal D(M),
\]
\[
D_{\lambda} \mathcal{K} \subseteq \mathcal{D}(M) \quad \text{ and } \quad   (\lambda-M)(-D_\lambda)f
  =
  f-f(\lambda)\mathbf1, \quad f \in \mathcal{K}.
\]
Consequently,
\[
  \ker(\lambda-M)=\{0\},
  \qquad
  \operatorname{Ran}(\lambda-M)=\ker\delta_\lambda,
  \qquad
  \operatorname{ind}(\lambda-M)=-1.
\]
Moreover,
\[
  (-D_\lambda)^*\mathcal K\subseteq\mathcal D(M^*),
  \qquad
  (\overline\lambda-M^*)(-D_\lambda)^*=I,
\]
and
\[
  \ker(\overline\lambda-M^*)
  =
  (\ker\delta_\lambda)^\perp,
  \qquad
  \operatorname{ind}(\overline\lambda-M^*)=1.
\]

\item
One has
\[
  \sigma(M)=\overline G,
  \qquad
  \rho(M)=\mathbb C\setminus\overline G,
\]
and, for \(\zeta\notin\overline G\), \(f \in \mathcal{K}\) and \(z \in G\),
\[
  ((\zeta-M)^{-1}f)(z)
  =
  \frac{f(z)}{\zeta-z},
  \qquad
  \|(\zeta-M)^{-1}\|
  \leq \frac{1}{\operatorname{dist}(\zeta,G)}.
\]
If \(G\) is bounded, then \(M\in\mathcal B(\mathcal K)\) and \(\|M\|\leq\sup_{z\in G}|z|\).
\end{enumerate}
\end{lemma}
Note that in the norm inequalities for \((\zeta-M)^{-1}\) and \(M\) one actually has equalities. However, this will not be needed below.

\begin{proof}
We adapt the weighted holomorphic \(L^2\)-space constructions of Mitiagin~\cite[Section~2.1]{Mityagin1978I} and den Boer~\cite[discussion after the proof of Theorem 1]{DenBoer1978} to the present setting. Put
\[
  w(z)
  =
  \frac{e^{-|z|^2}}
       {1+\sum_{j=1}^s|a_j(z)|^2},
  \qquad z\in G,
  \qquad
  d\mu=w\,dA,
\]
and let
\[
  \mathcal A:=H(G)\cap L^2(G,\mu)
\]
denote the corresponding weighted Bergman space, equipped with the \(L^2\)-inner product.

Since \(w\) is continuous and strictly positive, for every nonempty compact \(L\Subset G\) we may choose \(r>0\) such that
\[
  \{z:\operatorname{dist}(z,L)\le r\}\Subset G
\]
and put
\[
  m:=\min_{\operatorname{dist}(z,L)\le r}w(z)>0.
\]
Then, for every \(\lambda\in L\) and \(f\in\mathcal A\), the submean inequality applied to the subharmonic function \(|f|^2\) gives
\[
  |f(\lambda)|^2
  \le
  \frac{1}{\pi r^2}
  \int_{B(\lambda,r)}|f(z)|^2\,dA(z)
  \le
  \frac{1}{\pi r^2m}\|f\|_{L^2(\mu)}^2.
\]
Consequently,
\begin{equation}\label{eq:local-sup}
\max_{\lambda\in L}|f(\lambda)| \le \frac{1}{\sqrt{\pi r^2m}}\|f\|_{L^2(\mu)}, \qquad f\in\mathcal A.
\end{equation}
Hence convergence in \(\mathcal A\) implies locally uniform convergence
in \(G\).

To see that \(\mathcal A\) is closed in \(L^2(G,\mu)\), suppose
that \((f_n)\) is a sequence in \(\mathcal{A}\) that converges to some \(f\in L^2(G,\mu)\).
By the preceding estimate, \((f_n)\) is locally uniformly Cauchy
in \(G\) and therefore converges locally uniformly to a function
\(g\in H(G)\). On the other hand, a subsequence of \((f_n)\) converges
to \(f\) almost everywhere with respect to \(\mu\). Since this
subsequence also converges pointwise to \(g\), we have \(f=g\)
almost everywhere. Thus \(f\) has a representative in \(\mathcal A\),
showing that \(\mathcal A\) is closed in \(L^2(G,\mu)\) and hence
a Hilbert space.

Since \(w(z)\le e^{-|z|^2}\), we have
\[
  \mu(G)
  \le
  \int_{\mathbb C}e^{-|z|^2}\,dA(z)
  =
  \pi,
\]
so \(\mu\) is a finite Borel measure on the open set \(G \subseteq \mathbb C\) and the space \(L^2(G,\mu)\) is separable. Consequently, \(\mathcal A\) is separable as well.

Moreover, for every \(n\in\mathbb N_0\), we have
\[
  \int_G |z|^{2n}
  \left(1+\sum_{j=1}^s|a_j(z)|^2\right)\,d\mu(z)
  =
  \int_G |z|^{2n}e^{-|z|^2}\,dA(z)
  <\infty.
\]
Thus all polynomials and all functions \(z\mapsto z^na_j(z)\) belong to \(\mathcal A\).

Now define
\[
  \mathcal E
  =
  \{f\in\mathcal A:z^nf\in\mathcal A
    \text{ for every }n\in\mathbb N_0\}
\]
and its \(\mathcal{A}\)-closure
\[
  \mathcal K=\overline{\mathcal E}^{\,\mathcal A}.
\]
Then \(\mathcal K\) is separable, contains all polynomials and the functions \(a_1,\ldots,a_s\), and hence is infinite dimensional. It remains to verify the remaining assertions.

\textup{(i)} First of all, point evaluations on \(\mathcal K\) are bounded by \eqref{eq:local-sup}.

Now fix \(\lambda\in G\) and \(f\in\mathcal A\). To control \(D_\lambda f\), choose \(\varepsilon>0\) such that \(\overline{B(\lambda,\varepsilon)}\subseteq G\). Then, for \(|z-\lambda|\ge\varepsilon\),
\[
  |D_\lambda f(z)|
  \le
  \frac{1}{\varepsilon}
  \bigl(|f(z)|+|f(\lambda)|\bigr),
\]
whereas the maximum principle applied to the holomorphic function \(D_\lambda f\) gives
\[
  \max_{|z-\lambda|\le\varepsilon}|D_\lambda f(z)|
  \le
  \frac{2}{\varepsilon}
  \max_{|z-\lambda|\le\varepsilon}|f(z)|.
\]
Integrating these estimates separately over \(B(\lambda,\varepsilon)\) and its complement, and using \eqref{eq:local-sup} and \(\mu(G)<\infty\), we see that there exists \(C_\lambda>0\) such that
\[
  \|D_\lambda f\|_{L^2(\mu)}
  \le
  C_\lambda\|f\|_{L^2(\mu)},
  \qquad f\in\mathcal A.
\]
In particular, \(D_\lambda\in\mathcal B(\mathcal A)\). Moreover,
\[
  z^nD_\lambda f
  =
  D_\lambda(z^nf)-f(\lambda)D_\lambda(z^n).
\]
Since \(D_\lambda(z^n)\) is a polynomial, this identity shows that \(D_\lambda\mathcal E\subseteq\mathcal E\). By the continuity of \(D_{\lambda}\) on \(\mathcal{A}\) and the definition of \(\mathcal K\), it follows that \(D_\lambda\mathcal K\subseteq\mathcal K\). This concludes the proof of \textup{(i)}.

\textup{(ii)} Let \(M\) denote the maximal coordinate multiplication operator on \(\mathcal{K}\). Since \(z\mathcal E\subseteq\mathcal E\), we have \(\mathcal E\subseteq\mathcal D(M)\), so \(M\) is densely defined. If \((f_n)\) is a sequence in \(\mathcal D(M)\) and \(f,h\in\mathcal K\) are such that
\[
  f_n\to f,
  \qquad
  Mf_n\to h
\]
in \(\mathcal K\), both convergences are locally uniform by \eqref{eq:local-sup}, and hence
\[
  h(z)=zf(z),\qquad z\in G.
\]
Thus \(f\in\mathcal D(M)\) and \(Mf=h\), so \(M\) is closed. This proves \textup{(ii)}.

\textup{(iii)} For \(f\in\mathcal K\) and \(\lambda\in G\),
\[
  z(-D_\lambda f)
  =
  \lambda(-D_\lambda f)-f+f(\lambda)\mathbf1\in\mathcal K,
\]
and therefore \(D_\lambda\mathcal K\subseteq\mathcal D(M)\). Direct computation now gives
\[
  (\lambda-M)(-D_\lambda f)
  =
  f-f(\lambda)\mathbf1,
  \qquad f\in\mathcal K,
\]
and
\[
  (-D_\lambda)(\lambda-M)u=u,
  \qquad u\in\mathcal D(M).
\]
The second identity implies that \(\lambda-M\) is injective. The first identity shows that every function vanishing at \(\lambda\) belongs to the range of \(\lambda-M\), while every function in this range necessarily vanishes at \(\lambda\). Consequently,
\[
  \ker(\lambda-M)=\{0\},
  \qquad
  \operatorname{Ran}(\lambda-M)=\ker\delta_\lambda.
\]
Since \(\delta_\lambda\ne0\), the latter space is closed of codimension one, and therefore \(\operatorname{ind}(\lambda-M)=-1\).

Taking adjoints in \((-D_\lambda)(\lambda-M)=I\) gives
\[
  (\overline\lambda-M^*)(-D_\lambda)^*=I.
\]
Indeed, for \(u\in\mathcal D(M)\) and \(h\in\mathcal K\),
\[
  \langle(\lambda-M)u,(-D_\lambda)^*h\rangle
  =
  \langle u,h\rangle,
\]
which also verifies that \((-D_\lambda)^*h\in\mathcal D(M^*)\). Thus \(\overline\lambda-M^*\) is surjective, while
\[
  \ker(\overline\lambda-M^*)
  =
  \operatorname{Ran}(\lambda-M)^\perp
  =
  (\ker\delta_\lambda)^\perp
\]
is one-dimensional. Hence \(\operatorname{ind}(\overline\lambda-M^*)=1\), which concludes the proof of \textup{(iii)}.

\textup{(iv)} Now let \(\zeta\notin\overline G\) and let \(R_{\zeta}\) denote the maximal operator of multiplication by \((\zeta-z)^{-1}\) on \(\mathcal{A}\). Then  \(R_\zeta \in \mathcal{B}(\mathcal{A})\) and
\[
  \|R_\zeta\|
  \le
  \frac{1}{\operatorname{dist}(\zeta,G)}.
\]
Since \(z^nR_\zeta f=R_\zeta(z^nf)\), the operator \(R_{\zeta}\) preserves \(\mathcal E\), hence also \(\mathcal K\). Furthermore,
\[
  zR_\zeta f=\zeta R_\zeta f-f,
\]
so \(R_\zeta\mathcal K\subseteq\mathcal D(M)\), and
\[
  (\zeta-M)R_\zeta=I
  \quad\text{on }\mathcal K,
  \qquad
  R_\zeta(\zeta-M)=I
  \quad\text{on }\mathcal D(M).
\]
Thus
\[
  (\zeta-M)^{-1}=R_\zeta|_{\mathcal{K}},
  \qquad
  \|(\zeta-M)^{-1}\|
  \le
  \frac{1}{\operatorname{dist}(\zeta,G)}.
\]
Since \(G\subseteq\sigma(M)\) by \textup{(iii)} and the spectrum is closed,
\[
  \sigma(M)=\overline G,
  \qquad
  \rho(M)=\mathbb C\setminus\overline G.
\]
Finally, if \(G\) is bounded, put \(R:=\sup_{z\in G}|z|\).
Then multiplication by \(z\) is bounded on \(\mathcal A\) with
norm at most \(R\). Since \(z\mathcal E\subseteq\mathcal E\),
it follows that \(z\mathcal K\subseteq\mathcal K\).
Thus \(M\in\mathcal B(\mathcal K)\) and \(\|M\|\le R\).
This proves \textup{(iv)}.
\end{proof}

We continue with the proof of Proposition~\ref{prop:scalar-realization}. As already indicated in the introduction, we actually prove a stronger simultaneous realization result for finitely many holomorphic functions. It is this version that will be needed below.

\begin{theorem}\label{thm:simultaneous}
Let \(\Omega \subseteq \mathbb{C}\) be open and nonempty, let \(g_1,\ldots,g_n\in H(\Omega)\) and \(\beta_1,\ldots,\beta_n\in\mathbb C\), and put \(F=\mathbb C\setminus\Omega\). Then the following holds.
\begin{enumerate}
\item There are a separable infinite-dimensional Hilbert space \(\mathcal G\), a densely defined closed operator \(T\) on \(\mathcal G\), linearly independent vectors \(x_1,\ldots,x_n\in\mathcal G\), and a unit vector \(y\in\mathcal G\) such that
\[
  \rho(T)=\Omega,
  \qquad
  \sigma(T)=\sigma_{\mathrm{ess}}(T)=F,
\]
and, for \(j=1,\ldots,n\),
\[
\langle x_j,y\rangle=\beta_j \quad \text{ and } \quad   g_j(\lambda)
  =
  \langle(\lambda-T)^{-1}x_j,y\rangle, \quad \lambda \in \Omega.
\]

\item The operator \(T\) in \textup{(i)} can be chosen bounded if and only if
\(F\) is nonempty and compact and, on the unbounded component of
\(\Omega\),
\[
  g_j(\lambda)
  =
  \frac{\beta_j}{\lambda}+O(|\lambda|^{-2}),
  \qquad |\lambda|\to\infty,
\]
for \(j=1,\ldots,n\).
\end{enumerate}
\end{theorem}
In particular, the case \(n=1\) proves Proposition~\ref{prop:scalar-realization}.

\begin{remark}
  \label{rem:3}
  Note that the necessity of the conditions in \textup{(ii)} is immediate. Indeed,
  for a bounded \(T\) on an infinite-dimensional Hilbert space, the
  spectrum and essential spectrum must be compact and nonempty. Moreover,
  the resolvent expansion at infinity gives
  \[
    \langle(\lambda-T)^{-1}x_j,y\rangle
    =
    \frac{\langle x_j,y\rangle}{\lambda}
    +
    O(|\lambda|^{-2}),
    \qquad |\lambda|\to\infty.
  \]
\end{remark}
\begin{proof}
\textup{(i)} We use Lemma~\ref{lem:multiplication} to obtain Hilbert spaces \(\mathcal K := \mathcal K(g_1,\ldots,g_n,\Omega) \subseteq H(\Omega)\) and \(\mathcal K_\# := \mathcal K(\Omega^\#) \subseteq H(\Omega^{\#})\), where \(\Omega^{\#}= \{ \overline{\lambda} : \lambda \in \Omega \}\). The corresponding coordinate multiplication operators are denoted by \(M\) and \(M_\#\), respectively. Now set
\[
  e=\frac{\mathbf1}{\|\mathbf1\|}\in\mathcal K_\#
\]
and define the operator \(T_0\) on \(\mathcal K\oplus\mathcal K_\#\) by
\[
  T_0
  =
  \begin{pmatrix}
    M & \langle\,\cdot\,,e\rangle\mathbf1\\
    0 & M_\#^*
  \end{pmatrix},
  \qquad
  \mathcal D(T_0)
  =
  \mathcal D(M)\oplus\mathcal D(M_\#^*).
\]
Since \(M_\#\) is densely defined and closed, so is its adjoint \(M_\#^*\). Hence \(M\oplus M_\#^*\) is densely defined and closed by Lemma~\ref{lem:direct-sums}\textup{(ii)}. Since \(T_0\) differs from this operator only by the bounded rank-one operator in the upper right corner, \(T_0\) is densely defined and closed as well. We have already indicated at the beginning of this section why \(T_0\) is the correct choice. 

Let us first show that \(\Omega\subseteq\rho(T_0)\). To this end, fix \(\lambda\in\Omega\). Then by Lemma~\ref{lem:multiplication}\textup{(iii)},
\[
  S_\lambda:=-D_{\overline\lambda}^*
\]
is a bounded right inverse of \(\lambda-M_\#^*\), i.e.
\[
  (\lambda-M_\#^*)S_\lambda h=h, \quad h \in \mathcal{K}_{\#}.
\]
Moreover, since \(D_{\overline\lambda}e=0\), we have
\[
  \langle S_\lambda h,e\rangle=0,
  \qquad h\in\mathcal K_\#.
\]
Again by Lemma~\ref{lem:multiplication}\textup{(iii)}, \(\ker(\lambda-M_\#^*)\) is one-dimensional and
\[
  \ker(\lambda-M_\#^*)
  =
  \operatorname{Ran}(\overline\lambda-M_\#)^\perp
  =
  (\ker\delta_{\overline\lambda})^\perp.
\]
Since \(\delta_{\overline\lambda}(e)\ne0\), the vector \(e\) is not orthogonal to this kernel, and we may therefore choose \(v_\lambda\in\ker(\lambda-M_\#^*)\) such that
\[
  \langle v_\lambda,e\rangle=1.
\]
Consequently, for every \(h \in \mathcal{K}_{\#}\) the solutions \(v \in \mathcal{D}(M^{*}_{\#})\) of
\[
  (\lambda-M_\#^*)v=h
\]
are given by 
\begin{equation}\label{eq:2}
  v=S_\lambda h+cv_\lambda,
  \qquad c\in\mathbb C,
\end{equation}
and for such a solution \(\langle v,e\rangle=c\).

Now let \((f,h)\in\mathcal K\oplus\mathcal K_\#\) and \((u,v) \in \mathcal{D}(T_0)\). Then the equation
\[
  (\lambda-T_0)(u,v)=(f,h)
\]
is equivalent to
\[
  (\lambda-M_\#^*)v=h
\]
and
\[
  (\lambda-M)u=f+\langle v,e\rangle\mathbf1.
\]
Using the  description \eqref{eq:2} of the solutions \(v\) of the first equation, the second equation becomes
\[
  (\lambda-M)u=f+c\mathbf1.
\]
Since \(\operatorname{Ran}(\lambda-M)=\ker\delta_\lambda\), this equation is solvable if and only if
\[
  f(\lambda)+c=0.
\]
Hence necessarily \(c=-f(\lambda)\), so the two equations
determine \(v\) uniquely as
\[
  v=S_\lambda h-f(\lambda)v_\lambda.
\]
Moreover, the second equation transforms into
\[
  (\lambda-M)u=f-f(\lambda)\mathbf1,
\]
which, by Lemma~\ref{lem:multiplication}\textup{(iii)}, has the unique solution \(u=-D_\lambda f\). We have therefore obtained that \(\lambda-T_0 : \mathcal{D}(T_0) \to \mathcal{K} \oplus \mathcal{K}_{\#}\) is bijective with
\begin{equation}\label{eq:T0res}
(\lambda-T_0)^{-1}(f,h) = \bigl( -D_\lambda f,\, S_\lambda h-f(\lambda)v_\lambda \bigr).
\end{equation}
Since all operators and functionals  on the right are bounded, we obtain \(\Omega\subseteq\rho(T_0)\).

Next, we can provide a realization of \(g_1, \ldots, g_{n} \in \mathcal{K}\) in terms of \(T_0\). To this end, recall that \(\mathcal K_\#\), containing all polynomials, is infinite dimensional. We may therefore choose orthonormal vectors
\[
  w_1,\ldots,w_n\in e^\perp \subseteq \mathcal{K}_{\#}.
\]
Set
\[
  x_j^0=(g_j,-\beta_je+w_j),
  \qquad
  y_0=(0,-e),
  \qquad
  j=1,\ldots,n.
\]
Then \(\|y_0\|=1\) and
\[
  \langle x_j^0,y_0\rangle=\beta_j.
\]
The vectors \(x_1^0,\ldots,x_n^0\) are linearly independent, since the projection of a linear relation among them onto \(e^\perp\) gives a linear relation among \(w_1,\ldots,w_n\).

Using \eqref{eq:T0res}, we obtain, for \(\lambda\in\Omega\),
\[
  \begin{aligned}
  \langle(\lambda-T_0)^{-1}x_j^0,y_0\rangle
  =
  \left\langle
    S_\lambda(-\beta_je+w_j)
    -
    g_j(\lambda)v_\lambda,
    -e
  \right\rangle=
  g_j(\lambda),
  \end{aligned}
\]
because \(\langle S_\lambda h,e\rangle=0\) for all \(h\in\mathcal K_\#\) and \(\langle v_\lambda,e\rangle=1\).

It remains to prescribe the spectrum outside \(\Omega\). Suppose first that \(F\ne\emptyset\).

Choose a countable dense sequence \((\zeta_k)\) in \(F\), repeating each of its elements infinitely often, and define the normal diagonal \emph{padding operator} \(N_F\) on \(\ell^2(\mathbb N)\) by
\[
  N_F(a_k)=(\zeta_ka_k),
  \qquad
  \mathcal D(N_F)
  =
  \left\{(a_k)\in \ell^2: (\zeta_{k}a_{k}) \in \ell^2(\mathbb{N})\right\}.
\]
Then
\[
  \sigma(N_F)=\sigma_{\mathrm{ess}}(N_F)=F,
\]
and, for \(\lambda\notin F\),
\[
  \|(\lambda-N_F)^{-1}\|
  =
  \frac{1}{\operatorname{dist}(\lambda,F)}.
\]

Now put
\[
  \mathcal G=(\mathcal K\oplus\mathcal K_{\#})
  \oplus\ell^2(\mathbb N)
\]
and define
\[
  T=T_0\oplus N_F,
  \qquad y=y_0\oplus0,
  \qquad x_j=x_j^0\oplus0,
  \quad j=1,\ldots,n,
\]
with \(\mathcal D(T)=\mathcal D(T_0)\oplus\mathcal D(N_F)\).
By Lemma~\ref{lem:direct-sums}\textup{(i), (ii)}, \(\mathcal G\) is a separable infinite-dimensional Hilbert space,
and \(T\) is closed and densely defined.
Clearly, \(\|y\|=\|y_0\|=1\) and
\(x_1,\ldots,x_n\) are linearly independent.  Moreover, for every \(\lambda\in\Omega\), both diagonal blocks are invertible, so Lemma~\ref{lem:direct-sums}\textup{(iv)} gives \(\lambda\in\rho(T)\). On the other hand, if \(\lambda\in F\), then \(\lambda-N_F\) is not Fredholm, and therefore neither is \(\lambda-T\) by Lemma~\ref{lem:direct-sums}\textup{(v)}. Thus
\[
  \rho(T)=\Omega,
  \qquad
  \sigma(T)=\sigma_{\mathrm{ess}}(T)=F.
\]
Moreover, for every \(\lambda\in\Omega\) and \(j=1,\ldots,n\),
\[
  (\lambda-T)^{-1}x_j
  =
  (\lambda-T_0)^{-1}x_j^0\oplus0.
\]
Consequently,
\[
  \langle(\lambda-T)^{-1}x_j,y\rangle
  =
  \langle(\lambda-T_0)^{-1}x_j^0,y_0\rangle = g_{j}(\lambda),
  \qquad
  \langle x_j,y\rangle=\langle x_j^0,y_0\rangle = \beta_{j}.
\]

Finally, if \(F=\emptyset\), then \(\Omega=\mathbb C\), and we simply take \(\mathcal{G}= \mathcal{K} \oplus \mathcal{K}_{\#}\) and
\[
  T=T_0,\qquad y=y_{0}, \qquad x_j=x_j^0, \quad j=1, \ldots,n.
\]
Since \(\Omega\subseteq\rho(T_0)\), in this case
\[
  \rho(T)=\mathbb C,
  \qquad
  \sigma(T)=\sigma_{\mathrm{ess}}(T)=\emptyset.
\]
This proves the first part of the theorem.

\medskip

\noindent\textup{(ii)} The necessity of the conditions in \textup{(ii)} was
already noted in Remark~\ref{rem:3}. We now prove sufficiency. Assume therefore that \(F\) is nonempty and compact and that for \(j=1, \ldots,n\)
\[
  g_j(\lambda)
  =
  \frac{\beta_j}{\lambda}+O(|\lambda|^{-2}),
  \qquad |\lambda|\to\infty.
\]
We fix \(c\in F\) and consider the biholomorphic map
\[
  \phi:\mathbb C\setminus\{c\}
  \longrightarrow
  \mathbb C\setminus\{0\},
  \qquad
  \phi(\lambda)=(\lambda-c)^{-1}.
\]
Put
\[
  V:=\phi(\Omega)\cup\{0\}.
\]
Since \(F\) is compact, there exists \(R>0\) such that
\[
  \{\lambda\in\mathbb C:|\lambda|>R\}\subseteq\Omega.
\]
Hence \(\phi(\Omega)\) contains a punctured neighborhood of \(0\).
Thus \(0\) is an interior point of \(V\), and since
\(\phi(\Omega)\) is open, \(V\) is an open subset of \(\mathbb C\).

Now for \(\mu\in V\setminus\{0\}\) and \(j=1,\ldots,n\) we define
\[
  h_j(\mu)
  :=
  \frac{\beta_j\mu-g_j(c+\mu^{-1})}{\mu^2}.
\]
We claim that the singularity at \(0\) is removable. Indeed, as
\(\mu\to0\),
\[
  \frac{1}{c+\mu^{-1}}
  =
  \frac{\mu}{1+c\mu}
  =
  \mu+O(|\mu|^2),
\]
and therefore the assumed asymptotics of \(g_j\) give
\[
  g_j(c+\mu^{-1})
  =
  \beta_j\mu+O(|\mu|^2).
\]
Consequently,
\[
  h_j(\mu)=O(1),
  \qquad \mu\to0.
\]
By the removable singularity theorem, each \(h_j\) extends
holomorphically to \(0\), and we henceforth regard
\[
  h_j\in H(V),
  \qquad j=1,\ldots,n.
\]

We may now apply part \textup{(i)} to the open set \(V\), the
functions \(h_1,\ldots,h_n\), and the prescribed numbers
\(\beta_1,\ldots,\beta_n\). Thus there exist a separable
infinite-dimensional Hilbert space \(\mathcal G_0\), a densely
defined closed operator \(S\) on \(\mathcal G_0\), linearly
independent vectors \(x_1,\ldots,x_n\in\mathcal G_0\), and a unit
vector \(y\in\mathcal G_0\) such that
\[
  \rho(S)=V,
  \qquad
  \sigma(S)=\sigma_{\mathrm{ess}}(S)=\mathbb C\setminus V,
\]
and
\[
  h_j(\mu)
  =
  \langle(\mu-S)^{-1}x_j,y\rangle,
  \qquad
  \langle x_j,y\rangle=\beta_j,
  \qquad
  \mu\in V.
\]

Since \(0\in V=\rho(S)\), the inverse
\[
  Q:=S^{-1}
\]
is a bounded operator on \(\mathcal G_0\). We  define
\[
  T_0:=c+Q \in \mathcal{B}( \mathcal{G}_0).
\]
Now note that by spectral mapping, we have
\[
  \sigma(Q)\setminus\{0\}
  =
  \{\mu^{-1}:\mu\in\sigma(S)\},
\]
and the analogous identity holds for the essential spectrum. Since
\[
  \sigma(S)=\sigma_{\mathrm{ess}}(S)
  =
  \mathbb C\setminus V
  =
  \phi(F\setminus\{c\}),
\]
we thus obtain
\[
  \sigma(T_0)\setminus\{c\}
  =
  \sigma_{\mathrm{ess}}(T_0)\setminus\{c\}
  =
  F\setminus\{c\}.
\]
We will deal with the point \(c\) below. First, let us verify the required resolvent coefficients. To this end, let
\[
\lambda=\phi^{-1}(\mu)=\mu^{-1}+c,
\quad \mu\in V\setminus\{0\}=\phi(\Omega).
\]
From
\[
  \lambda-T_0
  =
  \mu^{-1}(I-\mu Q)
\]
and
\[
  (I-\mu Q)^{-1}
  =
  I-\mu(\mu-S)^{-1},
  \qquad \mu\in V\setminus\{0\},
\]
we obtain
\[
  (\lambda-T_0)^{-1}
  =
  \mu I-\mu^2(\mu-S)^{-1}.
\]
Hence, for \(\lambda\in\Omega\), and \(j=1, \ldots, n\)
\[
  \begin{aligned}
  \langle(\lambda-T_0)^{-1}x_j,y\rangle
  &=
  \mu\langle x_j,y\rangle
  -
  \mu^2
  \langle(\mu-S)^{-1}x_j,y\rangle
  \\
  &=
  \beta_j\mu-\mu^2h_j(\mu)
  \\
  &=
  g_j(c+\mu^{-1})
  =
  g_j(\lambda).
  \end{aligned}
\]

Now, to ensure that
\(c\) belongs to the essential spectrum as required, let
\[
  \mathcal G
  =
  \mathcal G_0\oplus\ell^2(\mathbb N)
\]
and set
\[
  T
  :=
  T_0\oplus cI_{\ell^2(\mathbb N)}.
\]
Then Lemma~\ref{lem:direct-sums}\textup{(iii), (iv), (v)} and the preceding spectral
relations imply that \(T \in \mathcal{B}(\mathcal{G})\) and
\[
  \rho(T)=\Omega,
  \qquad
  \sigma(T)=\sigma_{\mathrm{ess}}(T)=F.
\]

Finally, for \(j=1, \ldots, n\) put
\[
  \widehat x_j=x_j\oplus0,
  \qquad
  \widehat y=y\oplus0.
\]
Then \(\widehat x_1,\ldots,\widehat x_n\) are linearly independent,
\(\|\widehat y\|=1\), and
\[
  \langle\widehat x_j,\widehat y\rangle
  =
  \beta_j, \qquad  j=1, \ldots, n.
\]
Moreover, for every \(\lambda\in\Omega\),
\[
  \langle(\lambda-T)^{-1}\widehat x_j,\widehat y\rangle
  =
  \langle(\lambda-T_0)^{-1}x_j,y\rangle
  =
  g_j(\lambda), \qquad  j=1, \ldots,n.
\]
Thus \(T\) is a bounded realization with all the required
properties.
\end{proof}

\section{Proof of our main results}\label{sec:pair}

We first construct an operator that carries a prescribed index function \(\kappa:\Omega\to\mathbb Z\). As noted above, in the bounded one-operator setting, the corresponding realization problem is classical. See Chevreau~\cite[Theorem~2.3]{Chevreau1980} and, in particular, Conway~\cite[Proposition~5]{Conway1985}, as well as Remark~\ref{rem:spectral-picture}. We do not use these results in the proof below. Instead, we retain the direct construction from Section~\ref{sec:ingredients}, which treats arbitrary closed essential spectra and the unbounded case in the same framework.

\begin{proposition}\label{prop:index}
Let \(\Omega\subseteq\mathbb C\) be nonempty and open,
put \(F=\mathbb C\setminus\Omega\), and let
\(\kappa:\Omega\to\mathbb Z\) be locally constant.
Then there exist a separable Hilbert space \(\mathcal H_\kappa\) and a densely defined
closed operator \(R_\kappa\) on \(\mathcal H_\kappa\)
such that \(\sigma_{\mathrm{ess}}(R_\kappa)=F\),
\[
  \operatorname{ind}(\lambda-R_\kappa)=\kappa(\lambda),
  \qquad \lambda\in\Omega,
\]
and
\[
  \rho(R_\kappa)
  =\{\lambda\in\Omega:\kappa(\lambda)=0\}.
\]
If \(F\) is nonempty and compact and \(\kappa\) vanishes on
the unbounded component of \(\Omega\), then one may choose
\(R_\kappa\in\mathcal B(\mathcal H_\kappa)\).
\end{proposition}
\begin{proof}
Let \(\mathfrak C\) denote the at most countable set of connected
components of \(\Omega\). Since each \(U\in\mathfrak C\) is open and \(\kappa\) is locally constant, it is constant on \(U\). As above, we denote this constant
value by \(\kappa_U\). Our first aim is to construct an operator \(C_U\) such that \(\operatorname{ind}(\lambda-C_U)=\kappa_U\) for each \(\lambda \in U\).

To this end, for \(U\in\mathfrak C\) let \(M_U\) be the densely defined
closed multiplication operator on \(\mathcal K(U)\subseteq H(U)\)
provided by Lemma~\ref{lem:multiplication}. As in the proof of Theorem~\ref{thm:simultaneous}, we also
use the operator
\[
  P_U:=M_{U^\#}^*,
  \qquad
  U^\#:=\{\overline z:z\in U\},
\]
acting on the Hilbert space \(\mathcal K(U^\#)\).
From Lemma~\ref{lem:multiplication} we obtain
\[
  \operatorname{ind}(\lambda-M_U)=-1,
  \qquad
  \operatorname{ind}(\lambda-P_U)=1,
  \qquad \lambda\in U.
\]

Moreover, since
\[
  \operatorname{dist}(\overline\zeta,U^\#)
  =
  \operatorname{dist}(\zeta,U),
  \qquad
  \sup_{z\in U^\#}|z|
  =
  \sup_{z\in U}|z|,
\]
Lemma~\ref{lem:multiplication} also shows that, for
\(\zeta\notin\overline U\),
\[
  \zeta\in\rho(M_U)\cap\rho(P_U),
  \qquad
  \max\bigl\{
    \|(\zeta-M_U)^{-1}\|,
    \|(\zeta-P_U)^{-1}\|
  \bigr\}
  \le
  \frac{1}{\operatorname{dist}(\zeta,U)},
\]
and, if \(U\) is bounded, \(M_U\in\mathcal B(\mathcal K(U))\) and \(P_U\in\mathcal B(\mathcal K(U^\#))\), with
\[
  \max\bigl\{\|M_U\|,\|P_U\|\bigr\}
  \le
  \sup_{z\in U}|z|.
\]

Writing \(A^{\oplus n}\) for the \(n\)-fold direct sum of an operator
\(A\) with itself, we now define
\[
  C_U:=
  \begin{cases}
    P_U^{\oplus\kappa_U}, & \kappa_U>0,\\[1mm]
    M_U^{\oplus(-\kappa_U)}, & \kappa_U<0,\\[1mm]
    0_{\{0\}}, & \kappa_U=0,
  \end{cases}
\]
where \(0_{\{0\}}\) denotes the unique operator on the
zero-dimensional Hilbert space \(\{0\}\). In particular,
\[
  \rho(0_{\{0\}})=\mathbb C,
  \qquad
  \|(\lambda-0_{\{0\}})^{-1}\|=0,
  \qquad \lambda\in\mathbb C.
\]

By the preceding observations and
Lemma~\ref{lem:direct-sums}\textup{(i), (ii), (v)},
the operator \(C_U\) acts on a separable Hilbert space, is densely
defined and closed, and satisfies
\[
  \operatorname{ind}(\lambda-C_U)=\kappa_U,
  \qquad \lambda\in U.
\]
In particular,
\[
  \lambda \in \rho(C_U) \quad \Leftrightarrow \quad \kappa_U=0,
  \qquad \lambda\in U.
\]
Furthermore, Lemma~\ref{lem:direct-sums}\textup{(iii), (iv)} gives,
for \(\zeta\notin\overline U\),
\[
  \zeta\in\rho(C_U),
  \qquad
  \|(\zeta-C_U)^{-1}\|
  \le
  \frac{1}{\operatorname{dist}(\zeta,U)},
\]
and, if \(U\) is bounded,
\begin{equation}\label{eq:5}
  \|C_U\|
  \le
  \sup_{z\in U}|z|.
\end{equation}

Since the direct sum of all these \(C_U\) need not have the required essential spectrum, we add a padding operator. To this end, if \(F\ne\emptyset\), we choose \(N_F\) to be the padding operator on \(\ell^2(\mathbb{N})\) from the proof of Theorem~\ref{thm:simultaneous}. In particular, in this case
\[
  \sigma(N_F)=\sigma_{\mathrm{ess}}(N_F)=F,
  \qquad
  \|(\lambda-N_F)^{-1}\|
  =
  \frac{1}{\operatorname{dist}(\lambda,F)},
  \qquad \lambda\in\Omega.
\]
In addition, if \(F= \emptyset\), we choose \(N_F = N_{\emptyset} = 0_{\{0\}}\) to be the unique operator on the Hilbert space
\(\{0\}\).

We now put
\[
  R_\kappa
  :=
  N_F\oplus\bigoplus_{U\in\mathfrak C}C_U.
\]
Since there are at most countably many summands and each acts on a
separable Hilbert space, the underlying Hilbert space \(\mathcal H_\kappa\) is separable by Lemma~\ref{lem:direct-sums}\textup{(i)}. Moreover, by part \textup{(ii)} of the same lemma,
\(R_\kappa\) is densely defined and closed. It remains to determine the Fredholm and spectral properties of this operator.

To this end, fix \(\lambda\in\Omega\) and let \(U_\lambda \in \mathfrak C\) be the component of
\(\Omega\) containing \(\lambda\). Moreover, let
\begin{equation}\label{eq:4}
  S_\lambda
  :=
  N_F\oplus
  \bigoplus_{V\in\mathfrak C\setminus\{U_\lambda\}} C_V.
\end{equation}
To show that \(\lambda-S_{\lambda}\) is boundedly invertible, choose \(r>0\) such that
\[
  B(\lambda,r)\subseteq\Omega.
\]
Then the operator \(\lambda-N_F\) is boundedly invertible with \(\left\| (\lambda-N_F)^{-1}\right\| \leq r^{-1}\). Moreover, for each \(V \in \mathfrak C\setminus\{U_\lambda\}\) we have
\[
    \operatorname{dist}(\lambda,V)\ge r,
\]
since the ball \(B(\lambda,r)\) is connected and contained in \(U_{\lambda}\). Consequently, \(\lambda \notin \overline{V}\) and so \(\lambda-C_{V}\) is boundedly invertible with
\[
  \left\| (\lambda -C_V)^{-1}\right\| \leq \frac{1}{\operatorname{dist}(\lambda,V)} \leq \frac{1}{r}.
\]
Hence Lemma~\ref{lem:direct-sums}\textup{(iv)} implies that \(\lambda-S_\lambda\) is boundedly invertible.

Now we note that, up to the order of the summands, we have
\[
  \lambda-R_\kappa=(\lambda-C_{U_\lambda})\oplus (\lambda-S_\lambda).
\]
Hence Lemma~\ref{lem:direct-sums}\textup{(iv)} and \textup{(v)} imply that \(\lambda-R_{\kappa}\) is Fredholm with
\[
  \operatorname{ind}(\lambda-R_{\kappa})= \operatorname{ind}(\lambda-C_{U_{\lambda}})= \kappa_{U_{\lambda}} = \kappa(\lambda),
\]
and that
\[
  \lambda \in \rho(R_{\kappa}) \quad \Leftrightarrow \quad \lambda \in \rho(C_{U_{\lambda}}).
\]
Since the latter holds if and only if \(\kappa_{U_{\lambda}}=0\), we obtain
\[
  \lambda\in\rho(R_\kappa)
  \quad\Longleftrightarrow\quad
  \kappa(\lambda)=0.
\]

It remains to show that \(\sigma_{\mathrm {ess}}(R_{\kappa})=F\), i.e. \(\Phi(R_{\kappa})=\Omega\), and that the assertion concerning boundedness of \(R_{\kappa}\) is valid.

To this end, let us first note that we already proved above that \(\Omega \subseteq \Phi(R_{\kappa})\). In particular, if \(F=\emptyset\) we have equality. If \(F \neq\emptyset\) and \(\mu \in F\), then \(\mu-N_F\) is not a Fredholm operator, and so neither is \(\mu-R_{\kappa}\) by Lemma~\ref{lem:direct-sums}\textup{(v)}. Consequently, \(\Omega = \Phi(R_{\kappa})\) also in this case.

Finally, suppose that \(F\) is nonempty and compact and that \(\kappa\)
vanishes on the unbounded component \(U_\infty\) of \(\Omega\). We have to show that \(R_\kappa= N_F\oplus
  \bigoplus_{U\in\mathfrak C} C_U \) is bounded. To this end,
put
\[
  r_{0}:=\max_{z\in F}|z|.
\]
Then
\[
  \{z\in\mathbb C:|z|>r_{0}\} \subseteq U_{\infty},
\]
so every component \(U \in \mathfrak C \setminus \{ U_{\infty} \}\) is contained in \(B(0,r_0)\). In particular, from \eqref{eq:5} we obtain that
\[
  \forall U \in \mathfrak C \setminus \{ U_{\infty} \}: \quad \left\| C_U \right\| \leq r_0.
\]
Moreover,
\[
  \|N_F\|=\sup_{k\in\mathbb N}|\zeta_k|\le r_0,
  \qquad
  \|C_{U_\infty}\|=0.
\]
Hence, Lemma~\ref{lem:direct-sums}\textup{(iii)} implies that \(R_\kappa\in\mathcal B(\mathcal H_\kappa)\).
\end{proof}

We now pass from scalar resolvent coefficients to rank-one perturbation determinants. Let us briefly review the necessary material. For a more comprehensive account, see Kato~\cite[Chapter~IV, Section~6]{Kato}.

Let \(T\) be a densely defined closed operator on a Hilbert space \(\mathcal H\), and for \(x,y\in\mathcal H\) set
\[
  Kf=\langle f,y\rangle x,
  \qquad f\in\mathcal H.
\]
For \(\lambda\in\rho(T)\), we define the rank-one perturbation determinant of \(T\) by \(K\) as
\[
  d_{T,K}(\lambda)
  :=
  \det\bigl(I-K(\lambda-T)^{-1}\bigr).
\]
Since \(K(\lambda-T)^{-1}\) has rank at most one and \(K(\lambda-T)^{-1}f
  =
  \bigl\langle(\lambda-T)^{-1}f,y\bigr\rangle x\), the elementary rank-one determinant formula gives
\begin{equation}\label{eq:pert-det-scalar}
d_{T,K}(\lambda) = 1-\bigl\langle(\lambda-T)^{-1}x,y\bigr\rangle.
\end{equation}
In particular, \(d_{T,K}\) is holomorphic on \(\rho(T)\). Let us recall the standard relation between this determinant and the spectrum of \(T+K\).

\begin{proposition}\label{prop:det-spectrum}
Let \(T, K\) and \(d_{T,K}\) be as above.
\begin{enumerate}
\item
For every \(\lambda\in\rho(T)\),
\[
  \lambda\in\sigma(T+K)
  \quad\Longleftrightarrow\quad
  d_{T,K}(\lambda)=0.
\]

\item
If \(d_{T,K}(\lambda)=0\) for some \(\lambda \in \rho(T)\), then
\[
  \ker(\lambda-T-K)
  =
  \operatorname{span}\{(\lambda-T)^{-1}x\}.
\]
In particular, every eigenvalue of \(T+K\) lying in \(\rho(T)\) has geometric multiplicity one.

\item
Let \(U\) be a connected component of \(\rho(T)\) on which \(d_{T,K}\not\equiv0\). Then the zeros of \(d_{T,K}\) in \(U\) are precisely the discrete eigenvalues of \(T+K\) in \(U\), and the order of each zero equals the algebraic multiplicity of the corresponding eigenvalue.
\end{enumerate}
\end{proposition}

\begin{proof}
From the identity
\[
  \lambda-T-K
  =
  \bigl(I-K(\lambda-T)^{-1}\bigr)(\lambda-T),
  \qquad \lambda\in\rho(T),
\]
we see that \(\lambda\in\rho(T+K)\) if and only if
\(I-K(\lambda-T)^{-1}\) is invertible. Since
\(K(\lambda-T)^{-1}\) has rank at most one, this is equivalent to
\[
  d_{T,K}(\lambda) = \det\bigl(I-K(\lambda-T)^{-1}\bigr)\ne0,
\]
which proves \textup{(i)}.

For \textup{(ii)}, suppose that \(d_{T,K}(\lambda)=0\) and put
\[
  u=(\lambda-T)^{-1}x\in\mathcal D(T)=\mathcal D(T+K).
\]
Then \(\langle u,y\rangle=1\), so in particular \(u\ne0\). Moreover,
\[
  (\lambda-T-K)u
  =
  x-K(\lambda-T)^{-1}x
  =
  x-\bigl\langle(\lambda-T)^{-1}x,y\bigr\rangle x
  =
  d_{T,K}(\lambda)x
  =
  0.
\]
Conversely, if \(f\in\mathcal D(T)\) satisfies
\((\lambda-T-K)f=0\), then
\[
  f
  =
  (\lambda-T)^{-1}(\lambda-T)f
  =
  (\lambda-T)^{-1}Kf
  =
  \langle f,y\rangle(\lambda-T)^{-1}x
  \in
  \operatorname{span}\{u\}.
\]
Hence \(\ker(\lambda-T-K)
  =
  \operatorname{span}\{(\lambda-T)^{-1}x\}\) and the geometric multiplicity of every eigenvalue is one.

Part \textup{(iii)} follows from the Weinstein--Aronszajn formula;
see Kato~\cite[Chapter~IV, Section~6]{Kato}.
\end{proof}

We will also need the classical Weierstrass theorem for holomorphic functions; see, for instance, Rudin~\cite[Theorem~15.11]{Rudin}. In the form needed here, it states that for every nonempty open set \(U\subseteq\mathbb C\), every locally finite set \(\Lambda_U\subseteq U\), and every function \(m_U:\Lambda_U\to\mathbb N\), there exists a function \(h\in H(U)\) whose zeros are precisely the points of \(\Lambda_U\), and such that each \(\lambda\in\Lambda_U\) is a zero of order \(m_U(\lambda)\). Thus both the location and the multiplicity of the zeros can be prescribed arbitrarily, subject only to local finiteness.

We now use these observations in the proof of the following lemma, which
treats the case where the common index function \(\kappa\) of two given spectral
configurations vanishes on all of \(\Omega\).

\begin{lemma}\label{lem:zero-pair}
Let \(\Omega\subseteq\mathbb C\) be open and nonempty and let
\[
  \mathcal C_A=(\Omega_{A,\mathrm{full}},0,\Lambda_A,m_A),
  \qquad
  \mathcal C_B=(\Omega_{B,\mathrm{full}},0,\Lambda_B,m_B)
\]
be two zero-index spectral configurations. Then, for every
\(\alpha\in\mathbb C\), there exist a separable infinite-dimensional
Hilbert space \(\mathcal G\), a rank-one operator \(K\) on \(\mathcal G\)
with
\[
  \operatorname{tr}K=\alpha,
\]
and a densely defined closed operator \(A\) on \(\mathcal G\) such that
\(A\) realizes \(\mathcal C_A\) on \(\Omega\) and
\(B:=A+K\) realizes \(\mathcal C_B\) on \(\Omega\).

Moreover, if \(F=\mathbb C\setminus\Omega\) is nonempty and compact and
the set
\[
  \Omega_{A,\mathrm{full}}\cup\Lambda_A
  \cup\Omega_{B,\mathrm{full}}\cup\Lambda_B
\]
is bounded, then \(A\) may be chosen in \(\mathcal B(\mathcal G)\)
(and then \(B=A+K\in\mathcal B(\mathcal G)\) as well).
\end{lemma}

\begin{proof}
We first encode the two prescribed spectra by holomorphic functions.
To this end, define \(q_A\in H(\Omega)\) componentwise as follows.
If a component \(U\) of \(\Omega\) is not contained in
\(\Omega_{A,\mathrm{full}}\), then, by the Weierstrass theorem, there
exists a holomorphic function on \(U\) whose zeros are precisely the
points of \(\Lambda_A\cap U\), with orders prescribed by \(m_A\); choose
\(q_A|_U\) to be such a function. If
\(U\subseteq\Omega_{A,\mathrm{full}}\), put \(q_A|_U=0\).
Define \(q_B\) analogously, using
\(\Omega_{B,\mathrm{full}}\), \(\Lambda_B\), and \(m_B\).

Thus \(q_A,q_B\in H(\Omega)\), and on every component on which they
are not identically zero, their zeros are exactly the prescribed
discrete eigenvalues, with orders equal to their prescribed algebraic
multiplicities.

Now we apply Theorem~\ref{thm:simultaneous} to the functions
\[
  g_A:=1-q_A,
  \qquad
  g_B:=1-q_B,
\]
with the prescribed pairings \(0\) and \(\alpha\). In this way, we obtain
a closed densely defined operator \(T\) on a separable infinite-dimensional
Hilbert space \(\mathcal G\), a unit vector \(y\in\mathcal G\), and linearly
independent vectors \(x_A,x_B\in\mathcal G\) such that
\[
  \rho(T)=\Omega,
  \qquad
  \sigma(T)=\sigma_{\mathrm{ess}}(T)=F,
\]
and
\[
  g_A(\lambda)=\langle(\lambda-T)^{-1}x_A,y\rangle,
  \qquad
  g_B(\lambda)=\langle(\lambda-T)^{-1}x_B,y\rangle,
\]
with
\[
  \langle x_A,y\rangle=0,
  \qquad
  \langle x_B,y\rangle=\alpha.
\]

Now define bounded rank-one operators \(K_A,K_B\) on \(\mathcal G\) by
\[
  K_Af=\langle f,y\rangle x_A,
  \qquad
  K_Bf=\langle f,y\rangle x_B,
  \qquad f\in\mathcal G,
\]
and put
\[
  K:=K_B-K_A,
  \qquad
  A:=T+K_A,
  \qquad
  \mathcal D(A):=\mathcal D(T).
\]
Then \(A\) is densely defined and closed. Moreover,
\[
  Kf=\langle f,y\rangle(x_B-x_A),
  \qquad f\in\mathcal G,
\]
so \(K\) has rank one, since \(x_A\) and \(x_B\) are linearly independent.
By construction,
\[
  \operatorname{tr}K
  =
  \langle x_B-x_A,y\rangle
  =
  \alpha.
\]
Finally,
\[
  B:=A+K=T+K_B.
\]

Since finite-rank perturbations preserve the Fredholm set and the
Fredholm index, both \(A\) and \(B\) have essential spectrum \(F\) and
Fredholm index zero throughout \(\Omega\). Furthermore,
\eqref{eq:pert-det-scalar} gives
\[
  d_{T,K_A}=q_A,
  \qquad
  d_{T,K_B}=q_B.
\]

Now consider a component \(U\) of \(\Omega\) which is not contained in
\(\Omega_{A,\mathrm{full}}\). Then \(q_A|_U\not\equiv0\), and its zeros,
with orders, are precisely the points of \(\Lambda_A\cap U\) with their
prescribed algebraic multiplicities. Hence
Proposition~\ref{prop:det-spectrum}\textup{(i), (iii)} shows that the
spectrum of \(A\) in \(U\) consists precisely of these discrete
eigenvalues, with the prescribed algebraic multiplicities.
If instead \(U\subseteq\Omega_{A,\mathrm{full}}\), then
\(q_A|_U\equiv0\), so Proposition~\ref{prop:det-spectrum}\textup{(i)}
gives \(U\subseteq\sigma(A)\). Thus \(A\) realizes \(\mathcal C_A\) on \(\Omega\). The same argument, with \(q_B\) in place of \(q_A\), shows that \(B\) realizes \(\mathcal C_B\) on \(\Omega\).

It remains to prove the bounded assertion. We therefore assume that
\(F=\mathbb C\setminus\Omega\) is nonempty and compact and that
\[
  \Omega_{A,\mathrm{full}}\cup\Lambda_A
  \cup\Omega_{B,\mathrm{full}}\cup\Lambda_B
\]
is bounded. Let \(U_\infty\) denote the unbounded component of
\(\Omega\). Since both \(\Omega_{A,\mathrm{full}}\) and
\(\Omega_{B,\mathrm{full}}\) are bounded, \(U_\infty\) is not
contained in either of these sets.

We now repeat the preceding construction, choosing the functions
\(q_A\) and \(q_B\) on every component of \(\Omega\) different from
\(U_\infty\) exactly as above. On \(U_\infty\), however, we have to
choose them more carefully in order to obtain the asymptotics
required in Theorem~\ref{thm:simultaneous}\textup{(ii)}.

As in the proof of that theorem, fix \(c\in F\) and consider
\[
  \phi:\mathbb C\setminus\{c\}
  \longrightarrow
  \mathbb C\setminus\{0\},
  \qquad
  \phi(\lambda)=(\lambda-c)^{-1}.
\]
Put
\[
  V:=\phi(U_\infty)\cup\{0\}.
\]
Since \(F\) is compact, \(U_\infty\) contains a neighborhood of
\(\infty\), and hence \(V\) is an open set containing \(0\).

For \(C=A,B\), set
\[
  \Lambda_C'
  :=
  \phi(\Lambda_C\cap U_\infty)
  =
  \{(\lambda-c)^{-1}:
      \lambda\in\Lambda_C\cap U_\infty\}.
\]
Since \(\Lambda_C\cap U_\infty\) is locally finite in \(U_\infty\)
and
\[
  \phi:U_\infty\longrightarrow V\setminus\{0\}
\]
is biholomorphic, \(\Lambda_C'\) is locally finite in
\(V\setminus\{0\}\). Moreover, since
\(\Lambda_A\cup\Lambda_B\) is bounded, \(0\) cannot be an
accumulation point of \(\Lambda_C'\). Thus \(\Lambda_C'\) is locally
finite in \(V\).

By the Weierstrass theorem, for \(C=A,B\) there exists
\(\psi_C\in H(V)\) whose zeros are precisely the points of
\(\Lambda_C'\), with the zero at \((\lambda-c)^{-1}\) having order
\(m_C(\lambda)\). Since \(0\notin\Lambda_C'\), we have
\[
  \psi_C(0)\ne0.
\]

Now define
\[
  \eta_A:=0,
  \qquad
  \eta_B:=-\alpha.
\]
For \(C=A,B\), choose \(a_C\in\mathbb C\) such that
\[
  e^{a_C}=\frac{1}{\psi_C(0)},
\]
set
\[
  b_C
  :=
  \eta_C-\frac{\psi_C'(0)}{\psi_C(0)},
\]
and define
\[
  \varphi_C(\mu)
  :=
  \psi_C(\mu)e^{a_C+b_C\mu},
  \qquad
  \mu\in V.
\]
Then \(\varphi_C\) has the same zeros as \(\psi_C\), with the same
orders, and
\[
  \varphi_C(0)=1,
  \qquad
  \varphi_C'(0)=\eta_C.
\]

We now define \(q_A\) and \(q_B\) on \(U_\infty\) by
\[
  q_C(\lambda)
  :=
  \varphi_C\bigl(\phi(\lambda)\bigr),
  \qquad
  \lambda\in U_\infty,
  \qquad C=A,B.
\]
Since \(\phi:U_\infty\to V\setminus\{0\}\) is biholomorphic,
\(q_C\) is holomorphic on \(U_\infty\) and has precisely the zeros
\(\Lambda_C\cap U_\infty\), with orders prescribed by \(m_C\).
Thus the resulting functions \(q_A,q_B\in H(\Omega)\) have all the
zero-set properties used in the first part of the proof.

Moreover,
\[
  \phi(\lambda)
  =
  (\lambda-c)^{-1}
  =
  \lambda^{-1}+O(|\lambda|^{-2}),
  \qquad |\lambda|\to\infty,
\]
and
\[
  \varphi_C(\mu)
  =
  1+\eta_C\mu+O(|\mu|^2),
  \qquad \mu\to0.
\]
Consequently,
\[
  q_A(\lambda)
  =
  1+O(|\lambda|^{-2}),
\]
and
\[
  q_B(\lambda)
  =
  1-\frac{\alpha}{\lambda}
  +O(|\lambda|^{-2}),
  \qquad
  |\lambda|\to\infty,\quad \lambda\in U_\infty.
\]
Hence, for
\[
  g_A:=1-q_A,
  \qquad
  g_B:=1-q_B,
\]
we obtain
\[
  g_A(\lambda)
  =
  O(|\lambda|^{-2}),
  \qquad
  g_B(\lambda)
  =
  \frac{\alpha}{\lambda}
  +O(|\lambda|^{-2}),
  \qquad
  |\lambda|\to\infty,\quad \lambda\in U_\infty.
\]

These are precisely the asymptotics required in
Theorem~\ref{thm:simultaneous}\textup{(ii)} for the prescribed
pairings \(0\) and \(\alpha\). Applying that theorem to these choices
of \(g_A\) and \(g_B\), we may therefore choose \(T\) bounded.
Defining \(K_A\), \(K_B\), \(K\), \(A\), and \(B\) exactly as in the
first part of the proof, we then have
\[
  A=T+K_A\in\mathcal B(\mathcal G),
  \qquad
  B=T+K_B\in\mathcal B(\mathcal G).
\]
The preceding trace and determinant arguments apply without change.
Thus \(K\) has rank one and trace \(\alpha\), while \(A\) and \(B\)
realize \(\mathcal C_A\) and \(\mathcal C_B\), respectively. This
proves the bounded assertion.
\end{proof}

The next  lemma will allow us to realize spectral configurations with a \emph{given} rank-one operator.

\begin{lemma}\label{lem:similarity}
For \(j=1,2\), let \(K_j\) be a rank-one operator on a separable infinite-di\-men\-sional Hilbert space \(\mathcal H_j\). Then there is a bounded isomorphism \(S:\mathcal H_1\to\mathcal H_2\) such that
\[
  SK_1S^{-1}=K_2
\]
if and only if \(\operatorname{tr}K_1=\operatorname{tr}K_2\).
\end{lemma}
\begin{proof}
Similarity preserves the trace of a finite-rank operator, proving necessity. Conversely, suppose that \(\operatorname{tr}K_1=\operatorname{tr}K_2\).

For \(j=1,2\), we choose a two-dimensional subspace \(E_j\subseteq\mathcal H_j\) containing \(\operatorname{Ran}K_j\) and \(\operatorname{Ran}K_j^*\). Then \(E_j\) is reducing for \(K_j\). Moreover,
\[
  E_j^\perp
  \subseteq(\operatorname{Ran}K_j^*)^\perp
  =\ker K_j.
\]
Consequently, \(K_j\) vanishes on \(E_j^\perp\), and we have
\[
  K_j=M_j\oplus0_{E_j^\perp},
  \qquad M_j:=K_j|_{E_j},
\]
where \(M_j\) is a rank-one operator on the two-dimensional space \(E_j\), with
\[
  \operatorname{tr}M_j=\operatorname{tr}K_j.
\]

Now recall that two complex \(2\times2\) matrices of rank one are similar if and only if they have the same trace. Indeed, the Jordan normal form of a \(2\times2\) rank-one matrix  \(M\) is 
\[
  \begin{pmatrix} \operatorname{tr} M &0\\0&0\end{pmatrix}
  \quad\text{if } \operatorname{tr} M \ne0,
  \qquad
  \begin{pmatrix}0&1\\0&0\end{pmatrix}
  \quad\text{if } \operatorname{tr} M=0.
\]

Hence, since \(M_1\) and \(M_2\) have the same trace, we may choose an isomorphism \(S_0:E_1\to E_2\) such that
\[
  S_0M_1=M_2S_0.
\]
As \(E_1^\perp\) and \(E_2^\perp\) are separable infinite-dimensional Hilbert spaces, we may also choose a unitary isomorphism \(U:E_1^\perp\to E_2^\perp\). Then
\[
  S:=S_0\oplus U:\mathcal H_1\to\mathcal H_2
\]
is a bounded isomorphism satisfying
\[
  SK_1
  =(S_0M_1)\oplus0
  =(M_2S_0)\oplus0
  =K_2S.
\]
Hence \(SK_1S^{-1}=K_2\), completing the proof.
\end{proof}

Finally, we are prepared for the proof of Theorem~\ref{thm:main}.

\begin{proof}[Proof of Theorem~\ref{thm:main}]
Put \(F=\mathbb C\setminus\Omega\).

If \(\Omega=\emptyset\), we choose on the prescribed Hilbert space \(\mathcal H\) a closed normal operator \(A\) with
\[
  \sigma(A)=\sigma_{\mathrm{ess}}(A)=\mathbb C.
\]
For instance, after fixing an orthonormal basis, we may take a diagonal operator whose diagonal entries form a dense sequence in \(\mathbb C\), each entry repeated infinitely often. By Fredholm stability,
\[
  \sigma_{\mathrm{ess}}(A+K)
  =\sigma_{\mathrm{ess}}(A)=\mathbb C.
\]
Since there is no further spectral data to realize, this proves the theorem in case \(\Omega=\emptyset\).

Assume henceforth that \(\Omega\ne\emptyset\), and put
\[
  \Omega_0=\{\lambda\in\Omega:\kappa(\lambda)=0\}.
\]
Since every component of \(\Omega\) on which \(\kappa\ne0\) is contained in \(\Omega_{A, \mathrm{full}}\) and in \(\Omega_{B, \mathrm{full}}\), we have \(\Lambda_A\cup\Lambda_B\subseteq\Omega_0\). As a preliminary step, we now define two zero-index configurations on all of \(\Omega\) by setting
\[
  \mathcal C_A^0
    =\bigl(\Omega_{A,\mathrm{full}}\cap\Omega_0,
             0,\Lambda_A,m_A\bigr), \qquad  
  \mathcal C_B^0
    =\bigl(\Omega_{B,\mathrm{full}}\cap\Omega_0,
             0,\Lambda_B,m_B\bigr).
\]
By Lemma~\ref{lem:zero-pair}, applied with \(\alpha=\operatorname{tr}K\), there exist closed densely defined operators \(A_0\) and \(B_0\) on a separable infinite-dimensional Hilbert space such that
\[
  B_0=A_0+K_0,\qquad
  \operatorname{rank}K_0=1,\qquad
  \operatorname{tr}K_0=\operatorname{tr}K,
\]
and such that \(A_0\) realizes \(\mathcal{C}_{A}^{0}\) on \(\Omega\) and \(B_0\) realizes \(\mathcal{C}_{B}^{0}\) on \(\Omega\). In particular,
\[
    \sigma_{\mathrm{ess}}(A_0)
    =\sigma_{\mathrm{ess}}(B_0)=F
\]
and
\[
  \sigma(A_0)\cap\Omega
    =(\Omega_{A,\mathrm{full}}\cap\Omega_0)
      \cup\Lambda_A, \qquad 
  \sigma(B_0)\cap\Omega
    =(\Omega_{B,\mathrm{full}}\cap\Omega_0)
    \cup\Lambda_B.
  \]
Since \(\Lambda_{A} \cup \Lambda_{B} \subseteq \Omega_0\) it also follows that
\begin{equation}\label{eq:zero-part-resolvent}
  \Omega\setminus\Omega_0
  \subseteq\rho(A_0)\cap\rho(B_0).
\end{equation}

Next, we take the index carrier \(R_\kappa\) from
Proposition~\ref{prop:index}. Recall that this operator is densely
defined and closed on a separable Hilbert space and satisfies
\[
  \sigma_{\mathrm{ess}}(R_\kappa)=F,
  \qquad
  \operatorname{ind}(\lambda-R_\kappa)=\kappa(\lambda),
  \quad \lambda\in\Omega.
\]
Moreover, for every \(\lambda \in \Omega\)
\begin{equation}\label{eq:index-carrier-resolvent}
  \lambda\in\rho(R_\kappa)
  \quad\Longleftrightarrow\quad
  \kappa(\lambda)=0.
\end{equation}
Now let us set
\[
  \widehat A=R_\kappa\oplus A_0,
  \qquad
  \widehat B=R_\kappa\oplus B_0.
\]
By Lemma~\ref{lem:direct-sums}\textup{(i)} and \textup{(ii)}, both operators act on an infinite-dimensional separable Hilbert space and are densely defined and closed. Moreover,
\begin{equation}\label{eq:hat-perturbation}
  \widehat B=\widehat A+(0\oplus K_0).
\end{equation}

Let us verify that \(\widehat A\) realizes \(\mathcal C_A\).
First, by Lemma~\ref{lem:direct-sums}\textup{(v)},
\[
  \sigma_{\mathrm{ess}}(\widehat A)
  =\sigma_{\mathrm{ess}}(R_\kappa)
   \cup\sigma_{\mathrm{ess}}(A_0)
  =F.
\]
Moreover, since \(\lambda-A_0\) has Fredholm index zero for
\(\lambda \in \Omega\), the same lemma gives
\[
  \operatorname{ind}(\lambda-\widehat A)
  =\operatorname{ind}(\lambda-R_\kappa) =\kappa(\lambda),
  \qquad \lambda\in\Omega.
\]
Furthermore, from \eqref{eq:index-carrier-resolvent} and
Lemma~\ref{lem:direct-sums}\textup{(iv)} we obtain that, on the one hand
\[
  \sigma(\widehat A)\cap\Omega_0
  =\sigma(A_0)\cap\Omega_0
  =(\Omega_{A,\mathrm{full}}\cap\Omega_0)\cup\Lambda_A,
\]
and on the other hand \(\Omega \setminus \Omega_0 \subseteq \sigma(\widehat A)\). Since \(\Omega\setminus\Omega_0\subseteq\Omega_{A,\mathrm{full}}\),
we conclude that
\[
  \sigma(\widehat A)\cap\Omega
  =\Omega_{A,\mathrm{full}}\cup\Lambda_A.
\]

In particular, every \(\lambda\in\Lambda_A\) is an isolated
point of \(\sigma(\widehat A)\). Let us show that these points are discrete eigenvalues by considering the corresponding Riesz projections, see Definition \ref{def:riesz-projection}. For each \(\lambda \in \Lambda_{A}\) we have \(\kappa(\lambda)=0\), so
\(\lambda\in\rho(R_\kappa)\) and Lemma~\ref{lem:direct-sums}\textup{(vi)} gives
\[
  P_\lambda(\widehat A)=0\oplus P_\lambda(A_0),
  \qquad
  \operatorname{rank}P_\lambda(\widehat A)
  =\operatorname{rank}P_\lambda(A_0)
  =m_A(\lambda).
\]
Thus these points are discrete eigenvalues with the
prescribed algebraic multiplicities. This concludes the proof that \(\widehat A\)
realizes \(\mathcal C_A\) on \(\Omega\). The same argument, with \(B_0\) in place of \(A_0\), shows
that \(\widehat B\) realizes \(\mathcal C_B\) on \(\Omega\) as well.

It remains to replace the perturbation in \eqref{eq:hat-perturbation} by the prescribed operator \(K\). We observe that \(0\oplus K_0\) is a  rank-one operator with
\[
  \operatorname{tr}(0\oplus K_0)
  =\operatorname{tr}K_0
  =\operatorname{tr}K.
\]
Since the underlying Hilbert space is separable and infinite dimensional, Lemma~\ref{lem:similarity} yields a bounded isomorphism \(S\) onto the prescribed Hilbert space \(\mathcal H\) such that
\[
  S(0\oplus K_0)S^{-1}=K.
\]
We set
\[
  A=S\widehat A S^{-1},
  \qquad
  B=S\widehat B S^{-1}.
\]
Then \(A\) and \(B\) are closed densely defined operators, and
\[
  B-A
  =S(\widehat B-\widehat A)S^{-1}
  =S(0\oplus K_0)S^{-1}
  =K.
\]
Since bounded similarity preserves the spectrum, the essential spectrum, Fredholm indices, and algebraic multiplicities of isolated eigenvalues, we conclude that \(A\) and \(B=A+K\) realize the prescribed spectral configurations.
\end{proof}

Next we consider the proof of the existence of bounded realizations.

\begin{proof}[Proof of Corollary~\ref{cor:bounded-main}]
Suppose first that \(A\) and \(A+K\) are bounded realizations on an infinite-dimensional Hilbert space. It is well known that in this case their spectra are compact and their joint essential spectrum \(F\) is compact and nonempty. Hence
\[
  F\cup\Omega_{A,\mathrm{full}}\cup\Lambda_A
  \cup\Omega_{B,\mathrm{full}}\cup\Lambda_B
  =F\cup\sigma(A)\cup\sigma(B)
\]
is bounded.

Conversely, assume that \(F\ne\emptyset\) and that the set on the left is bounded. Choose \(R>0\) such that this set is contained in \(B(0,R)\). In particular, \(F\) is compact and
\[
  \Omega_{A,\mathrm{full}}\cup\Omega_{B,\mathrm{full}}
  \subseteq B(0,R).
\]
Denote the unbounded component of \(\Omega=\mathbb C\setminus F\) by \(U_\infty\). Since \(\Omega_{A,\mathrm{full}}\) and \(\Omega_{B,\mathrm{full}}\) are unions of bounded components of \(\Omega\), we have
\[
  U_\infty\cap
  (\Omega_{A,\mathrm{full}}\cup\Omega_{B,\mathrm{full}})
  =\emptyset.
\]
Since
\[
  \{\lambda\in\Omega:\kappa(\lambda)\ne0\}
  \subseteq
  \Omega_{A,\mathrm{full}}\cap\Omega_{B,\mathrm{full}},
\]
we conclude that
\[
  \kappa|_{U_\infty}=0.
\]
Thus the hypotheses of the bounded parts of both Lemma~\ref{lem:zero-pair} and Proposition~\ref{prop:index} are satisfied: the prescribed spectral data are bounded, \(F\) is compact and nonempty, and the common index vanishes on the unbounded component.

We may therefore carry out the construction in the proof of Theorem~\ref{thm:main} using the bounded zero-index realizations from Lemma~\ref{lem:zero-pair} and the bounded index carrier from Proposition~\ref{prop:index}. By Lemma~\ref{lem:direct-sums}\textup{(iii)}, their direct sums give bounded operators \(\widehat A\) and \(\widehat B\) realizing the prescribed configurations before the final similarity. Since this similarity is implemented by a boundedly invertible operator, it preserves boundedness and produces the required bounded realizations \(A\) and \(A+K\).
\end{proof}

We now provide the proof of Corollary~\ref{cor:finite-rank}, the extension to a prescribed \(K\) of finite rank.

\begin{proof}[Proof of Corollary~\ref{cor:finite-rank}]
We use strong induction on \(n=\operatorname{rank}K\), the case
\(n=1\) following from Theorem~\ref{thm:main} and
Corollary~\ref{cor:bounded-main}. Suppose that \(n\ge2\) and that
the conclusion holds for any given Fredholm-compatible pair of
spectral configurations on \(\Omega\) and every prescribed
perturbation of positive rank less than \(n\) on a given
infinite-dimensional separable Hilbert space.

Now let
\[
  \mathcal C_{A}=(\Omega_{A,\mathrm{full}},\kappa,\Lambda_{A},m_{A}),
  \qquad
  \mathcal C_{B}=(\Omega_{B,\mathrm{full}},\kappa,\Lambda_{B},m_{B})
\]
be two spectral configurations on \(\Omega\) with the same index
function, and let \(K\) have rank \(n\) on a prescribed separable
infinite-dimensional Hilbert space \(\mathcal H\).

First, let us choose a finite-dimensional \(K\)-invariant subspace
\(W\) of \(\mathcal{H}\) such that
\[
  \operatorname{rank}(K|_W)=1.
\]
If \(K\) has a nonzero eigenvalue, we take \(W\) to be the span
of a corresponding eigenvector. Otherwise the restriction of
\(K\) to the finite-dimensional invariant space
\(\operatorname{Ran}K\) has only the eigenvalue zero and is
therefore nilpotent.
Since \(K\mathcal H\subseteq\operatorname{Ran}K\), it follows
that \(K\) itself is nilpotent. We can thus find \(r\ge2\) such that
\(K^r=0\) and \(K^{r-1}\ne0\). Now we choose any \(u\in\mathcal H\)
with \(K^{r-1}u\ne0\) and set
\[
  W=\operatorname{span}\{K^{r-2}u,K^{r-1}u\}.
\]
We now use \(W\) to obtain an upper triangular block
representation of \(K\), with one diagonal block of rank one
and the other of rank at most \(n-1\). To apply the induction
hypothesis, both blocks must act on infinite-dimensional
spaces. We therefore enlarge \(W\) by a suitable subspace
of \(\ker K\), keeping its orthogonal complement infinite
dimensional.

To this end, let us first note that, since \(K\) has finite
rank and \(W^{\perp}\) is infinite dimensional, the space
\(\ker K\cap W^\perp = \ker(K|_{W^\perp})\) is infinite
dimensional as well.
Splitting an orthonormal basis into two infinite parts we can
write this subspace as an orthogonal sum
\[
  \ker K\cap W^\perp=H_0\oplus H_1
\]
with closed infinite-dimensional subspaces \(H_0,H_{1}\).

Now let us set
\[
  \mathcal M:=W\oplus H_0,
  \qquad
  \mathcal N:=\mathcal M^\perp\supseteq H_1.
\]

Then both spaces are separable and infinite dimensional.
Since \(W\) is \(K\)-invariant and \(K\) vanishes on \(H_0\),
we have \(K\mathcal M=KW\subseteq\mathcal M\). Thus, with
respect to the decomposition
\(\mathcal H=\mathcal M\oplus\mathcal N\), we obtain
\[
  K=
  \begin{pmatrix}
    L&C\\
    0&K_0
  \end{pmatrix},
  \qquad
  L:=K|_{\mathcal M},
  \qquad
  \operatorname{rank}L=1.
\]
In particular, note that \(K_{0}\in\mathcal{B}(\mathcal{N})\)
and \(\operatorname{rank}K_0<n\).

Applying the induction hypothesis to the rank-one operator
\(L\), we now obtain a closed densely defined operator \(A_1\)
on \(\mathcal M\) such that \(A_1\) and \(B_1:=A_1+L\)
realize \(\mathcal C_A\) and \(\mathcal C_B\), respectively.

In addition, on \(\mathcal N\), we choose a closed densely
defined operator \(D\) such that \(D\) and \(D+K_0\) realize
the identical trivial spectral configurations
\((\emptyset,0,\emptyset,\emptyset)\) on \(\Omega\).
For \(K_0\ne0\), this is possible by the induction hypothesis.
For \(K_0=0\), we take a single such operator \(D\), whose
existence follows from Theorem~\ref{thm:main}, applied with
an arbitrary rank-one perturbation on \(\mathcal N\).
In the bounded case, we use
Corollary~\ref{cor:bounded-main} as well.
In either case,
\[
  \rho(D)=\rho(D+K_0)=\Omega,
  \qquad
  \sigma_{\mathrm{ess}}(D)=\mathbb{C}\setminus\Omega=:F.
\]

Now set
\[
  A=A_1\oplus D,
  \qquad
  B=A+K=
  \begin{pmatrix}
    B_1&C\\
    0&D+K_0
  \end{pmatrix}.
\]
By Lemma~\ref{lem:direct-sums}\textup{(ii)}, \(A\) is densely
defined and closed. The same holds for \(B\), since \(K\)
is bounded. We now verify that \(A\) and \(B\) realize the
prescribed spectral configurations.

To this end, let us first note that
Lemma~\ref{lem:direct-sums}\textup{(v)} gives
\[
  \sigma_{\mathrm{ess}}(A)
  =
  \sigma_{\mathrm{ess}}(A_1)
  \cup\sigma_{\mathrm{ess}}(D)
  =F,
\]
so Fredholm stability under finite-rank perturbations also yields
\[
  \sigma_{\mathrm{ess}}(B)=F.
\]
Moreover, since \(\lambda-D\) is invertible for every
\(\lambda\in\Omega\), Lemma~\ref{lem:direct-sums}\textup{(v)}
and invariance of the Fredholm index under finite-rank
perturbations also imply that
\[
  \operatorname{ind}(\lambda-B)
  =\operatorname{ind}(\lambda-A)
  =\operatorname{ind}(\lambda-A_1)
  =\kappa(\lambda),
  \qquad \lambda\in\Omega.
\]
Now let us consider the spectrum of \(A\) in \(\Omega\). Since \(\rho(D)=\Omega\),
Lemma~\ref{lem:direct-sums}\textup{(iv)} shows that
\[
  \sigma(A)\cap\Omega
  =\sigma(A_1)\cap\Omega
  =\Omega_{A,\mathrm{full}}\cup\Lambda_A.
\]
In particular, every \(\lambda\in\Lambda_A\) is an isolated
point of \(\sigma(A)\). Moreover,
Lemma~\ref{lem:direct-sums}\textup{(vi)} gives
\[
  P_\lambda(A)=P_\lambda(A_1)\oplus0,
  \qquad
  \operatorname{rank}P_\lambda(A)
  =\operatorname{rank}P_\lambda(A_1)=m_A(\lambda),
  \quad \lambda\in\Lambda_A.
\]
Thus every \(\lambda\in\Lambda_A\) is a discrete eigenvalue
of \(A\) with the prescribed algebraic multiplicity. This concludes the proof that \(A\) realizes \(\mathcal{C}_A\)
on \(\Omega\).

Now let us turn to \(B\). To this end, for \(\lambda\in\Omega\) we set
\(Q_\lambda=\lambda-D-K_0\). Then, by construction,
\(Q_\lambda\) is boundedly invertible, and
\[
  \lambda-B
  =
  \begin{pmatrix}
    I&-CQ_\lambda^{-1}\\
    0&I
  \end{pmatrix}
  \begin{pmatrix}
    \lambda-B_1&0\\
    0&Q_\lambda
  \end{pmatrix}.
\]
The first factor is boundedly invertible, its inverse being
obtained by changing the sign of the upper-right block.
Consequently,
\[
  \sigma(B)\cap\Omega
  =\sigma(B_1)\cap\Omega
  =\Omega_{B,\mathrm{full}}\cup\Lambda_B.
\]
In particular, each \(\lambda\in\Lambda_B\) is an isolated
point of \(\sigma(B)\). It remains to show that it is actually
a discrete eigenvalue of \(B\) with algebraic multiplicity
\(m_{B}(\lambda)\). To this end, let us choose a small positively
oriented circle \(\Gamma\) around \(\lambda\), whose closed
interior is contained in \(\Omega\) and intersects
\(\sigma(B_1)\) only in \(\lambda\).
For \(z\in\Gamma\), the above factorization gives
\[
  (z-B)^{-1}
  =
  \begin{pmatrix}
    (z-B_1)^{-1}&(z-B_1)^{-1}C Q_z^{-1}\\
    0&Q_z^{-1}
  \end{pmatrix}.
\]
Since \(Q_z^{-1}\) is holomorphic on \(\Omega\), its integral
over \(\Gamma\) vanishes. Integrating the resolvent formula
and using Definition~\ref{def:riesz-projection}, we thus obtain
\[
  P_\lambda(B)=
  \begin{pmatrix}
    P_\lambda(B_1)&X\\
    0&0
  \end{pmatrix},
  \qquad
  X:=\frac{1}{2\pi i}\int_\Gamma
       (z-B_1)^{-1}C Q_z^{-1}\,dz.
\]
Here \(X:\mathcal N\to\mathcal M\) is bounded, since the
integrand is continuous in operator norm on \(\Gamma\).
As \(P_\lambda(B)^2=P_\lambda(B)\), we have
\(P_\lambda(B_1)X=X\). Consequently,
\[
  \operatorname{Ran}P_\lambda(B)
  =
  \operatorname{Ran}P_\lambda(B_1)\oplus\{0\},
\]
and hence
\[
  \operatorname{rank}P_\lambda(B)
  =
  \operatorname{rank}P_\lambda(B_1)
  =m_B(\lambda).
\]
Thus \(\lambda\) is a discrete eigenvalue of \(B\) with
\(m_a(\lambda;B)=m_B(\lambda)\). This concludes the proof
that \(B\) realizes \(\mathcal{C}_B\) on \(\Omega\).

Finally we note that, under the conditions of
Corollary~\ref{cor:bounded-main}, we can use bounded
realizations throughout the induction.
Then \(A_1\) and \(D\) are bounded, and Lemma~\ref{lem:direct-sums}\textup{(iii)} shows that \(A=A_1\oplus D\) is bounded. Since \(B=A+K\), the same holds for \(B\).
Conversely, bounded realizations necessarily have bounded
spectra and nonempty common essential spectrum,
independently of the rank of \(K\).
\end{proof}

We conclude this article with the announced corollary concerning prescribed eigenvalue asymptotics.

\begin{corollary}\label{cor:rate}
Let \(F\subsetneq\mathbb C\) be nonempty and closed, let \(\mathcal H\) be a separable infinite-di\-men\-sional Hilbert space, let \(K\in\mathcal B(\mathcal H)\) have rank one and let \((t_n)_{n \in \mathbb{N}}\) be a strictly decreasing sequence of positive numbers with \(t_n\to0\). Then there exists a densely defined closed operator \(A\) on \(\mathcal H\) with \(\sigma(A)=\sigma_{\mathrm {ess}}(A)=F\), such that \(\sigma_{d}(A+K)\) consists of algebraically simple eigenvalues \(\lambda_n\), \(n\in\mathbb N\), and 
\[
  \operatorname{dist}(\lambda_n,F)=t_n
  \qquad\text{for all sufficiently large }n.
\]
Moreover, if \(F\) is compact, \(A\) may be chosen bounded and the distance equality holds for every \(n \in \mathbb{N}\).
\end{corollary}
Note that, in particular, for \(t_n=1/\log(n+1)\) we obtain \eqref{eq:1}.
\begin{proof}
We choose \(\xi\notin F\) and \(\zeta\in F\) with
\[
  \delta:=\operatorname{dist}(\xi,F)=|\xi-\zeta|>0.
\]
If \(F\) is compact, we may choose \(\xi\) with \(\delta>t_1\). Choose \(N \in \mathbb{N}\) such that \(t_n<\delta\) for \(n\ge N\), taking \(N=1\) in the compact case, and extend \((t_n)_{n\ge N}\) to a strictly decreasing sequence \((s_n)_{n\ge1}\) in \((0,\delta)\). Set
\[
  \lambda_n:=\zeta+\frac{s_n}{\delta}(\xi-\zeta),
  \qquad n\ge1.
\]
Note that \(|\lambda_n-\zeta|=s_n\) and \(|\lambda_n-\xi|=\delta-s_n\), since \(\lambda_n\) lies on the segment between \(\zeta\) and \(\xi\). Hence, using the triangle inequality, we obtain
\[
  s_n=\delta-(\delta-s_n)
  =\operatorname{dist}(\xi,F)-|\xi-\lambda_n|
  \le \operatorname{dist}(\lambda_n,F)
  \le |\lambda_n-\zeta|=s_n.
\]
Thus \(\operatorname{dist}(\lambda_n,F)=s_n\), which equals \(t_n\) for \(n\ge N\). Moreover, the \(\lambda_n\) are pairwise distinct, and \(\Lambda=\{\lambda_n:n\ge1\}\) is locally finite in \(\Omega=\mathbb C\setminus F\), since \(s_n\to0\).

Applying Theorem~\ref{thm:main} on \(\Omega\) with vanishing index, no full components, no discrete spectrum prescribed for \(A\), and with \(\Lambda\) and multiplicity one prescribed for \(A+K\) gives the assertions in the general case. If \(F\) is compact, then \(N=1\) and \(s_n=t_n\) for every \(n\). Since \(\Lambda\) is bounded, Corollary~\ref{cor:bounded-main} allows \(A\) to be chosen bounded. 
\end{proof}

\appendix
\section{Direct sums of closed operators}\label{app:direct-sums}

This appendix collects some facts about direct sums of closed operators. While all of them are standard, we were unable to find a suitable reference in the literature.

Since it will be needed below as well as elsewhere in the article, we also recall the following definition of Riesz projections.

\begin{definition}\label{def:riesz-projection}
Let \(T\) be a closed operator on a Hilbert space \(\mathcal{H}\), and suppose that
\[
\sigma(T)=\sigma\mathbin{\dot\cup}\tau,
\]
where \(\sigma\) is nonempty and contained in a bounded open set \(\Delta \subseteq \mathbb{C}\) whose boundary consists of finitely many piecewise smooth Jordan curves and satisfies \(\overline{\Delta}\cap\tau=\emptyset\). We call \(\sigma\) an isolated spectral subset of \(T\) and define its Riesz projection by
\[
P_\sigma(T):=\frac{1}{2\pi i}\int_{\partial\Delta}(z-T)^{-1}\,dz,
\]
where the boundary is positively oriented; see \cite[Chapter~XV]{GohbergGoldbergKaashoek1990}. This indeed defines a projection on \(\mathcal{H}\) and the definition is independent of the choice of \(\Delta\). We also set \(P_\emptyset(T):=0\).

For an isolated point \(\lambda\in\sigma(T)\), we write \(P_\lambda(T):=P_{\{\lambda\}}(T)\). This projection is nonzero, and if it has finite rank, then \(\lambda\) is called a discrete eigenvalue of \(T\) with algebraic multiplicity
\[
m_a(\lambda;T):=\operatorname{rank}P_\lambda(T).
\]
\end{definition}

\begin{lemma}\label{lem:direct-sums}
Let \((\mathcal H_j)_{j\in J}\) be a finite or countable nonempty family of Hilbert spaces, and for each \(j\in J\) let \(T_j\) be a closed densely defined operator on \(\mathcal H_j\).

We set
\[
\mathcal H:=\bigoplus_{j\in J}\mathcal H_j,
\qquad
\langle x,y\rangle_{\mathcal H}
:=\sum_{j\in J}\langle x_j,y_j\rangle_{\mathcal H_j},
\]
and define the direct sum operator \(T:=\bigoplus_{j\in J}T_j\) by
\[
Tx=(T_jx_j)_{j\in J},\qquad
\mathcal D(T)=
\left\{
x=(x_j)_{j\in J}\in\mathcal H:
x_j\in\mathcal D(T_j),\
\sum_{j\in J}\|T_jx_j\|^2<\infty
\right\}.
\]
Then the following hold.

\begin{enumerate}
\item \(\mathcal H\) is a Hilbert space. If all \(\mathcal H_j\) are separable, then so is \(\mathcal H\).

\item \(T\) is densely defined and closed.

\item \(T\in\mathcal B(\mathcal H)\) if and only if \(T_j\in\mathcal B(\mathcal H_j)\) for every \(j\in J\) and
\[
\sup_{j\in J}\|T_j\|<\infty.
\]
In that case, \(\|T\|=\sup_{j\in J}\|T_j\|\).

\item
\[
\lambda\in\rho(T)
\quad\Longleftrightarrow\quad
\lambda\in \bigcap_{j \in J} \rho(T_j) 
\quad\text{and}\quad
\sup_{j\in J}\|(\lambda-T_j)^{-1}\|<\infty.
\]
In that case,
\[
(\lambda-T)^{-1}
=\bigoplus_{j\in J}(\lambda-T_j)^{-1}.
\]

\item The following are equivalent:
\begin{enumerate}
\item \(T\) is Fredholm.
\item \(T_j\) is Fredholm for every \(j\in J\), and there is a finite set \(J_0\subseteq J\) such that \(T_j\) is boundedly invertible for \(j\in J\setminus J_0\) and the set \(\{\|T_j^{-1}\|:j\in J\setminus J_0\}\) is bounded.
\end{enumerate}
In that case,
\[
\operatorname{ind}T
=\sum_{j\in J_0}\operatorname{ind}T_j
=\sum_{j\in J}\operatorname{ind}T_j.
\]

\item Let \(\sigma\) be an isolated spectral subset of \(T\), and put \(\sigma_j:=\sigma\cap\sigma(T_j)\). Then
\[
P_\sigma(T)=\bigoplus_{j\in J}P_{\sigma_j}(T_j).
\]
Consequently,
\[
\operatorname{rank}P_\sigma(T)
=\sum_{j\in J}\operatorname{rank}P_{\sigma_j}(T_j),
\]
where both sides may be infinite.

In particular, if \(\lambda\) is an isolated point of \(\sigma(T)\), then \(\lambda\) is a discrete eigenvalue of \(T\) if and only if
\[
\sum_{j\in J}\operatorname{rank}P_\lambda(T_j)<\infty,
\]
and in that case,
\[
m_a(\lambda;T)=\sum_{j\in J}m_a(\lambda;T_j).
\]
Here we set \(P_\lambda(T_j)=0\) and \(m_a(\lambda;T_j)=0\) whenever \(\lambda\in\rho(T_j)\).
\end{enumerate}
\end{lemma}

\begin{remark}\label{rem:4}
If \(J\) is finite, the uniform resolvent bound in \textup{(iv)} is automatic. Moreover, \(T\) is Fredholm if and only if every \(T_j\) is Fredholm; indeed, in \textup{(v)} we may take \(J_0=J\).
\end{remark}

\begin{proof}
We denote by \(\iota_j:\mathcal H_j\to\mathcal H\) the canonical
isometric embedding. By definition, for \(u\in\mathcal D(T_j)\),
we have \(\iota_j u\in\mathcal D(T)\) and
\(T\iota_j u=\iota_jT_j u\).

\textup{(i)} The completeness of \(\mathcal H\) is standard; see Conway \cite[Chapter I, Section 6, Proposition 6.2]{Conway1990}. If all \(\mathcal H_j\) are separable, finite-support vectors with entries in fixed countable dense subsets form a countable dense subset of \(\mathcal H\).

\textup{(ii)} Finite-support vectors with entries in \(\mathcal D(T_j)\) are dense in \(\mathcal H\), so \(T\) is densely defined.

Now suppose that \((x^{(n)})_{n \in \mathbb{N}}\) is a sequence in \(\mathcal D(T)\) and that \(x^{(n)}\to x\)  and \(Tx^{(n)}\to y\) in \(\mathcal{H}\). Then
\[
x_j^{(n)}\to x_j,\qquad
T_jx_j^{(n)}\to y_j \quad \text{ in } \mathcal{H}_{j}
\]
for every \(j\in J\). Since \(T_j\) is closed, \(x_j\in\mathcal D(T_j)\) and \(T_jx_j=y_j\). As \(y\in\mathcal H\), it follows that \(x\in\mathcal D(T)\) and \(Tx=y\). Thus \(T\) is closed.

\textup{(iii)} We first suppose that every \(T_j\) is bounded and that \(C:=\sup_{j\in J}\|T_j\|<\infty\). Then
\[
\sum_{j\in J}\|T_jx_j\|^2
\le C^2\sum_{j\in J}\|x_j\|^2
=C^2\|x\|^2
\qquad(x\in\mathcal H).
\]
Hence \(\mathcal D(T)=\mathcal H\) and \(\|T\|\le C\). Conversely, suppose that \(T\in\mathcal B(\mathcal H)\). Then \(\mathcal D(T)=\mathcal H\), so \(\mathcal D(T_j)=\mathcal H_j\) and \(T_j=\iota_j^*T\iota_j\). Consequently, \(\|T_j\|\le\|T\|\) for every \(j\), which also proves the norm identity.

\textup{(iv)} We first suppose that \(\lambda\in\rho(T_j)\) for every \(j\) and
\[
C:=\sup_{j\in J}\|(\lambda-T_j)^{-1}\|<\infty.
\]
By \textup{(iii)}, the operator
\[
R:=\bigoplus_{j\in J}(\lambda-T_j)^{-1}
\]
is bounded on \(\mathcal H\). For \(y\in\mathcal H\), we set \(x=Ry\). Then
\[
T_jx_j=\lambda x_j-y_j,\qquad j\in J.
\]
Since \(\lambda x-y\in\mathcal H\), we have \(x\in\mathcal D(T)\) and \((\lambda-T)Ry=y\). Conversely, for \(x\in\mathcal D(T)\), the componentwise resolvent identities give \(R(\lambda-T)x=x\). Thus \(R=(\lambda-T)^{-1}\).

Now suppose that \(\lambda\in\rho(T)\). The identity
\[
(\lambda-T)\iota_j u=\iota_j(\lambda-T_j)u,
\qquad u\in\mathcal D(T_j),
\]
shows that \(\lambda-T_j\) is injective. For \(y\in\mathcal H_j\),
set \(x=(\lambda-T)^{-1}\iota_jy\). Then
\(x_j\in\mathcal D(T_j)\) and \((\lambda-T_j)x_j=y\), so
\(\lambda-T_j\) is also surjective. Moreover,
\[
\|(\lambda-T_j)^{-1}y\|
=\|x_j\|
\le\|x\|
\le\|(\lambda-T)^{-1}\|\,\|y\|.
\]
Thus \(\lambda\in\rho(T_j)\), and taking the supremum over \(j\)
gives the required uniform bound.

\textup{(v)} For arbitrary \(J\), the componentwise definition of \(T\) gives
\begin{equation}\label{eq:direct-sum-defects}
\ker T=\bigoplus_{j\in J}\ker T_j, \qquad (\operatorname{Ran}T)^\perp=\bigoplus_{j\in J}(\operatorname{Ran}T_j)^\perp.
\end{equation}

We first consider finite direct sums. In this case, we also have
\[
\operatorname{Ran}T
=\bigoplus_{j\in J}\operatorname{Ran}T_j.
\]
Since a finite orthogonal direct sum of subspaces is closed if and only if each summand is closed, the preceding identities show that \(T\) is Fredholm if and only if every \(T_j\) is Fredholm. In that case,
\[
\operatorname{ind}T
=\sum_{j\in J}\operatorname{ind}T_j.
\]
This proves the assertion for finite \(J\), taking \(J_0=J\).

Now we suppose that \(J\) is infinite. If a finite set \(J_0\)
with the stated properties exists, we set
\[
S:=\bigoplus_{j\in J_0}T_j,
\qquad
U:=\bigoplus_{j\in J\setminus J_0}T_j.
\]
Then \(S\) is Fredholm by the finite case, and \(U\) is boundedly
invertible by \textup{(iv)}. Since \(T=S\oplus U\), including the
corresponding decomposition of its domain, the finite case
shows that \(T\) is Fredholm. Here and below, empty summands
are omitted.

Conversely, suppose that \(T\) is Fredholm. For \(y\in\operatorname{Ran}T_j\), clearly \(\iota_jy\in\operatorname{Ran}T\). Conversely, if \(Tx=\iota_jy\), then \(T_jx_j=y\). Hence
\[
\iota_j(\operatorname{Ran}T_j)
=\operatorname{Ran}T\cap\iota_j(\mathcal H_j).
\]
Since \(\operatorname{Ran}T\) is closed, so is each \(\operatorname{Ran}T_j\). The decompositions in \eqref{eq:direct-sum-defects} therefore imply that every \(T_j\) is Fredholm and that only finitely many summands have a nontrivial kernel or cokernel. We let \(J_0\) consist of the corresponding indices. Then for \(j\notin J_0\), the operator \(T_j\) is bijective and, being closed, has a bounded inverse. It remains to show the uniform boundedness of these inverses.

Since \(T\) and its range are closed, the inverse of
\[
T:\mathcal D(T)\cap(\ker T)^\perp
\longrightarrow\operatorname{Ran}T
\]
has closed graph as an operator from \(\operatorname{Ran}T\)
to \((\ker T)^\perp\). The closed graph theorem therefore
gives \(c>0\) such that
\[
\|Tx\|\ge c\|x\|,
\qquad
x\in\mathcal D(T)\cap(\ker T)^\perp.
\]
For \(j\notin J_0\) and \(y\in\mathcal H_j\), we set \(x=\iota_jT_j^{-1}y\). Since \(\ker T_j=\{0\}\), the kernel decomposition \eqref{eq:direct-sum-defects} shows that \(x\perp\ker T\). Thus
\[
\|y\|=\|Tx\|\ge c\|x\|
=c\|T_j^{-1}y\|,
\]
and consequently
\[
\sup_{j\in J\setminus J_0}\|T_j^{-1}\|
\le c^{-1}<\infty.
\]
Finally, using \eqref{eq:direct-sum-defects}, we obtain
\[
\operatorname{ind}T
=\sum_{j\in J}\operatorname{ind}T_j
=\sum_{j\in J_0}\operatorname{ind}T_j.
\]

\textup{(vi)} We choose \(\Delta\) as in Definition~\ref{def:riesz-projection}
and set \(\Gamma=\partial\Delta\). By \textup{(iv)}, \(\sigma(T_j)\subseteq\sigma(T)\) and \(\Gamma\subseteq\rho(T_j)\) for every \(j\in J\). Thus
\[
\sigma_j=\sigma(T_j)\cap\Delta=\sigma\cap\sigma(T_j)
\]
is either empty or an isolated spectral subset of \(T_j\).

Again using \textup{(iv)}, we have
\[
(z-T)^{-1}=\bigoplus_{j\in J}(z-T_j)^{-1},
\qquad z\in\Gamma,
\]
and
\[
\sup_{j\in J,\ z\in\Gamma}\|(z-T_j)^{-1}\|
\le\sup_{z\in\Gamma}\|(z-T)^{-1}\|<\infty.
\]
Integrating coordinatewise over \(\Gamma\), we obtain the stated
projection formula. Taking ranks then gives the rank identity,
including the case of infinite rank. Finally, for \(\sigma=\{\lambda\}\), we obtain the remaining assertions from Definition~\ref{def:riesz-projection} and the rank identity.
\end{proof}

\hypertarget{ai-declaration}{}
\section*{Declaration on the use of generative AI}

The development of this paper began in early July 2026, when I asked ChatGPT 5.6 Sol to look through my list of open research questions and suggest possible approaches. Much to my surprise at the time, one of its suggestions was a substantial generalization of our result \cite{HansmannKatriel2013} that every Cauchy transform on the unit disk can be realized using rank-one perturbation determinants: namely, it argued that every holomorphic function on the disk can be realized in this way. In trying to understand and verify the suggested proof and trace the origins of its underlying ideas, I was led to the classical realization theory of Mitiagin~\cite{Mityagin1978I,Mityagin1978II} and den Boer~\cite{DenBoer1978}, of which I had previously been unaware. One question led to another, and the content of the paper gradually took shape through a series of interactions with ChatGPT (OpenAI) and occasional exchanges with Claude (Anthropic).

I independently verified all mathematical statements, proofs, and references appearing in this manuscript. Based on preliminary drafts generated by ChatGPT, the final version was written by me and reflects my own decisions and preferences. I take full responsibility for its mathematical content and exposition.

\end{document}